\documentclass[11pt]{article}
\usepackage{graphicx}
\usepackage{amssymb}
\usepackage{amsmath}
\usepackage{subfigure}
\usepackage{booktabs}
\usepackage{subfigure}

\title{Preconditioned smoothers for the full approximation scheme for the RANS equations}
\author{Philipp Birken$^{\mbox{\tiny\rm 1}}$, Jonathan Bull$^{\mbox{\tiny\rm 2}}$, Antony Jameson$^{\mbox{\tiny\rm 3}}$}
\begin{document}
\maketitle
\baselineskip=0.9
\normalbaselineskip
\vspace{-3pt}
\begin{center}{\footnotesize\em $^{\mbox{\tiny\rm 1}}$Centre for the
    mathematical sciences, Numerical Analysis, Lund University, Lund, Sweden\\ email: philipp.birken\symbol{'100}na.lu.se\\
$^{\mbox{\tiny\rm 2}}$Division of Scientific Computing, Dept. of Information Technology, Uppsala University, 
Box 337, 75105 Uppsala, Sweden\\
$^{\mbox{\tiny\rm 2}}$Stanford University, Department of Aeronautics \& Astronautics, Stanford, CA 94305, USA}
\end{center}

\begin{abstract}
We consider multigrid methods for finite volume discretizations of the Reynolds Averaged Navier-Stokes (RANS) equations for both steady and unsteady flows. We analyze the effect of different smoothers based on pseudo time iterations, such as explicit and additive Runge-Kutta (ARK) methods. Furthermore, we derive the new class of additive W (AW) methods from Rosenbrock smoothers. This gives rise to two classes of preconditioned smoothers, preconditioned ARK and additive W (AW), which are implemented the exact same way, but have different parameters and properties. The new derivation allows to choose some of these based on results for time integration methods. As preconditioners, we consider SGS preconditioners based on flux vector splitting discretizations with a cutoff function for small eigenvalues. We compare these methods based on a discrete Fourier analysis. Numerical results on pitching and plunging airfoils identify AW3 as the best smoother regarding overall efficiency. Specifically, for the NACA 64A010 airfoil steady-state convergence rates of as low as 0.85 were achieved, or a reduction of 6 orders of magnitude in approximately 25 pseudo-time iterations. Unsteady convergence rates of as low as 0.77 were achieved, or a reduction of 11 orders of magnitude in approximately 70 pseudo-time iterations.
\end{abstract}


\section{Introduction}

We are interested in numerical methods for compressible wall bounded turbulent flows as they appear in many problems in industry. Therefore, both steady and unsteady flows will be considered. Numerically, these are characterized by strong nonlinearities and a large number of unknowns, due to the requirement of resolving the boundary layer. High fidelity approaches such as Direct Numerical Simulation (DNS) or Large Eddy Simulation (LES) are slowly getting within reach through improvements in high order discretization methods. Nevertheless, these approaches are, and will remain in the foreseeable future, far too costly to be standard tools in industry. 

However, low fidelity turbulence modelling based on the Reynolds Averaged Navier-Stokes (RANS) equations discretized using second order finite volume discretizations is a good choice for many industrial problems where turbulence matters. For steady flows, this comes down to solving one nonlinear system. In the unsteady case, the time discretization has to be at least partially implicit, due to the extremely fine grids in the boundary layer, requiring solving one or more nonlinear systems per time step. The choice for numerical methods for these comes down to Jacobian-free Newton-Krylov (JFNK) methods with appropriate preconditioners or nonlinear multigrid methods (Full Approximation scheme - FAS) with appropriate smoothers, see \cite{birkenhabil} for an overview. 

In this article, we focus on improving the convergence rate of agglomeration multigrid methods, which are the standard in the aeronautical industry. For the type of problems considered here, two aspects have been identified that affect solver efficiency. Firstly, the flow is convection dominated. Secondly the grid has high aspect ratio cells. It is important to note that the viscous terms in the RANS equations do not pose problems in itself. Instead, it is that they cause the boundary layer to appear, thus making high aspect ratio grids necessary. These aspects are shared by the Euler equations, meaning that solvers developed for one equation may also be effective for the other.

With regards to convection dominated flows, smoothers such as Jacobi or Gau\ss -Seidel do not perform well, in particular when the flow is aligned with the grid \cite{mulder:89}. One idea has been to adjust multigrid restriction and prolongation by using directional or semi coarsening that respects the flow direction \cite{mulder:92}. This approach has the problem to be significantly more complicated to implement than standard agglomeration. Thus, the alternative is to adjust the smoother. As it turned out, symmetric Gau\ss -Seidel (SGS) is an excellent smoother for the Euler equations even for grid aligned flow \cite{caujam:01}, simply because it takes into account propagation of information in the flow direction and backwards. 

However, when discretizing the Euler equations on high aspect ratio grids suitable for wall bounded viscous flows, this smoother does not perform well. During the last ten years, the idea of preconditioned pseudo time iterations has garnered interest \cite{swturo:07, langer:11, lanli:11, langer:13, langer:14, lasckr:15, jameso:15, bibuja:16b, bibuja:16a}. This goes back to the additive Runge-Kutta (ARK) smoothers originally introduced in \cite{copsay:80} and independently in a multigrid setting in \cite{jameso:85}. These exhibit slow convergence, but if they are combined with a preconditioner, methods result that work well for high Reynolds number high aspect ratio RANS simulations. 

The preconditioned RK method suggested in \cite{swturo:07} was properly derived in \cite{langer:13}. There it is shown that this smoother arises from using a Rosenbrock method and then approximating the system matrix in the arising linear systems. This is in fact called a W method in the literature on ordinary differential equations. Consequently, we now introduce the class of preconditioned additive W methods and derive preconditioned additive explicit Runge-Kutta methods. This allows us to identify the roles the preconditioners have to play. As for preconditioners itself, it turns out that again, SGS is a very good choice, as reported in \cite{swturo:07,langer:13}. 

The specific convergence rate attainable depends on the discretization, in particular the flux function. Here, we consider the Jameson-Schmidt-Turkel (JST) scheme in its latest version \cite{jameso:15b}. We perform a discrete Fourier analysis of the smoother for the linearized Euler equations on grids with variable aspect ratios. This is justified, since the core issues of convection and high aspect ratio grids are present in this problem. 

A convenient truth is that if we have a fast steady-state solver then it can be used to build a fast unsteady solver via dual timestepping. However, there are subtle differences that affect convergence and stability. In particular, the eigenvalues of the amplification matrix are scaled and shifted in the unsteady case relative to the steady case. For a fuller discussion of these issues we refer to our earlier work \cite{bibuja:16b,bibuja:16a} and to \cite{bedods:17}. We compare the analytical behavior and numerical performance of iterative smoothers for steady and unsteady problems.

The article is structured as follows. We first present the governing equations and the discretization, then we describe multigrid methods and at length the smoothers considered. Then we present a Fourier analysis based on the Euler equations and finally numerical results for airfoil test cases. 

\section{Discretization}

We consider the two dimensional compressible (U)RANS equations, where the vector of conservative variables is $(\rho, \rho v_1, \rho v_2, \rho E)^T$ and the convective and viscous fluxes are given by

\begin{equation*}
    F_i^c = \left (
    \begin{tabular}{c}
    	$\rho v_i$ \\
    	$\rho v_i v_1 + p\delta_{i1}$ \\
    	$\rho v_i v_1 + p\delta_{i1}$ \\
    	$\rho v_i H$
    \end{tabular}
    \right ), \quad
    F_i^v = \left (
    \begin{tabular}{c}
    	0 \\
    	$\tau_{i1}$ \\
    	$\tau_{i2}$ \\
    	$v_j \tau_{ij} + \frac{\mu+\mu_t}{Pr} (C_p \partial_i T)$
    \end{tabular}
    \right ), \quad i=1,2,
\end{equation*}
\begin{equation*}
    \tau_{ij} = (\mu+\mu_t) (\partial_{x_j} v_i + \partial_{x_i} v_j - \tfrac 2 3 \delta_{ij} \partial_{x_k} v_k),
\end{equation*}
\begin{equation*}
    q_j = \left( \frac{\mu}{Pr}+\frac{\mu_t}{Pr_t} \right) \partial_{x_j} \left( H - \tfrac 1 2 v_k v_k \right)
\end{equation*}
where we used the Einstein notation. 

Here, $\rho$ is the density, $v_i$ the velocity components and $E$ the total energy per unit mass. The enthalpy is given by $H = E + p/\rho$ with $p = (\gamma -1)\rho(E-1/2 v_k v_k)$ being the pressure and $\gamma=1.4$ the adiabatic index for an ideal gas. Furthermore, $\tau_{ij}$ is the sum of the viscous and Reynolds stress tensors, $q_j$ the sum of the diffusive and turbulent heat fluxes, $\mu$ the dynamic viscosity, $\mu_t$ the turbulent viscosity and $Pr$, $Pr_t$ the dynamic Prandtl and turbulent Prandtl numbers. 

As a turbulence model, we use the 0-equation Baldwin-Lomax model \cite{ballom:78} for two reasons. Firstly, it performs well for flows around airfoils we use as primary motivation. Secondly, an algebraic turbulence model does not lead to additional questions regarding implementation as 1- or 2-equation models do. We believe that these difficulties have to be systematically looked at, but separately from this investigation.

The equations are discretized using a finite volume method on a structured mesh and the JST scheme as flux function. There are many variants of this method, see e.g. \cite{jameso:15b}. Here, we use the following, for simplicity written as if for a one dimensional problem:
\[{\bf f}_{j+1/2}^{JST}({\bf u}) = \frac12({\bf f}^R(\bar{{\bf u}}_j)+{\bf f}^R(\bar{{\bf u}}_{j+1})) + {\bf d}_{j+1/2}({\bf u}).\]
Here, ${\bf f}^R(\bar{{\bf u}})$ is the difference of the convective and the viscous fluxes, ${\bf u}\in\mathbb{R}^m$ is the vector of all discrete unknowns and $\bar{{\bf u}}=(\rho, \rho v, \rho E)$ is the vector of conservative variables. The artificial viscosity terms are given by
\[{\bf d}_{j+1/2}({\bf u})=\epsilon_{j+1/2}^{(2)}\Delta {\bf w}_j - \epsilon_{j+1/2}^{(4)}(\Delta {\bf w}_{j+1}-2\Delta {\bf w}_j+\Delta {\bf w}_{j-1})\]
with $\Delta_j$ being the forward difference operator and the vector ${\bf w}$ being $\bar{{\bf u}}$ where in the last component, the energy density has been replaced by the enthalpy density. 

The scalar coefficient functions $\epsilon^{(2)}_{j+1/2}$ and $\epsilon^{(4)}_{j+1/2}$ are given by
\begin{equation}\label{e2}
\epsilon_{j+1/2}^{(2)} = s_{j+1/2}r_{j+1/2}
\end{equation}
and
\begin{equation}\label{e4}
\epsilon^{(4)}_{j+1/2} = \max(0, r_{j+1/2}/32 - 2\epsilon_{j+1/2}^{(2)}).
\end{equation}
Here, the entropy sensor $s_{j+1/2}=\min(0.25,\max(s_j, s_{j+1}))$ given via
\[s_j=\left | \frac{S_{j+1}-2S_j + S_{j-1}}{S_{j+1}+2S_j + 2S_{j-1} +0.001}\right |\]
with $S = p/\rho^{\gamma}$. For the Euler equations, it is suggested to instead use a corresponding pressure sensor. 

Furthermore, $r_{i+1/2}$ is the scalar diffusion coefficient, which approximates the spectral radius and is chosen instead of a matrix valued diffusion as in other versions of this scheme. It is
\[r_{j+1/2} = \max(r_j,r_{j+1}).\] 
The specific choice of $r_j$ is important with respect to stability and the convergence speed of the multigrid method. Here, we use the locally largest eigenvalue $r_j = |v_{n_j}|+a_j$ as a basis, where $a$ is the speed of sound. In the multidimensional case, this is further modified to be \cite{martin:87}:
\begin{align}\label{spectralcorrection}
\tilde r_i = r_i(1+(r_j/r_i)^{2/3}), \\
\tilde r_j = r_j(1+(r_i/r_j)^{2/3}), \nonumber
\end{align} 
where $r_i$ corresponds to the $x$ direction and $r_j$ to the $y$ direction. 

Additionally, to obtain velocity and temperature gradients needed for the viscous fluxes, we exploit that we have a cell centered method on a structured grid and use dual grids around vertices to avoid checker board effects \cite[p. 40]{jameso:04}. 

For boundary conditions, we use the no slip condition at fixed wall and far field conditions at outer boundaries. These are implemented using Riemann invariants \cite[p. 38]{jameso:04}. 

In time, we use BDF-2 with a fixed time step $\Delta t$, resulting at time $t_{n+1}$ in an equation system of the form

\begin{equation}\label{nonlin-eq}
{\bf F}({\bf u}) := \frac{3{\bf u}-4{\bf u}^n + {\bf u}^{n-1}}{2\Delta t} + {\bf \Omega}^{-1}{\bf f}({\bf u}) = {\bf 0}.
\end{equation}
Here, ${\bf f}({\bf u})$ describes the spatial discretization, whereas ${\bf \Omega}$ is a diagonal matrix with the volumes of the mesh cells as entries. We thus obtain
\[\frac{\partial {\bf F}}{\partial {\bf u}} = \frac{3}{2\Delta t} {\bf I} + {\bf \Omega}^{-1} \frac{\partial {\bf f}}{\partial {\bf u}}.\]

For a steady state problem, we just have
\begin{equation}\label{nonlin-eq-steady}
{\bf F}({\bf u}) := {\bf \Omega}^{-1} {\bf f}({\bf u}) = {\bf 0}.
\end{equation}

\section{The full approximation scheme}

As mentioned in the introducion, we use an agglomeration FAS to solve equations \eqref{nonlin-eq} and \eqref{nonlin-eq-steady}. To employ a multigrid method, we need a hierarchical sequence of grids with the coarsest grid being denoted by level $l=0$. The coarse grids are obtained by agglomerating 4 neighboring cells to one. On the coarse grids, the problem is discretized using a first order version of the JST scheme that does not use fourth order differences or an entropy sensor. 

The iteration is performed as a W-cycle, where on the coarsest grid, one smoothing step is performed. This gives the following pseudo code:

Function FAS-W-cycle$({\bf u}_l, {\bf s}_l, l)$
\begin{itemize}
\item ${\bf u}_l = {\bf S}_l^{\nu_1}({\bf u}_l,{\bf s}_l)$ (Presmoothing)
\item if $(l>0)$
\begin{itemize}
\item ${\bf r}_l={\bf s}_l-{\bf F}_l({\bf u}_l)$
\item ${\bf \tilde{u}}_{l-1} = {\bf R}_{l-1,l} {\bf u}_l$ (Restriction of solution)
\item ${\bf s}_{l-1}={\bf F}_{l-1}({\bf \tilde{u}}_{l-1}) + {\bf R}_{l-1,l}{\bf r}_l$ (Restriction of residual)
\item For $j=1,2$: call FAS-W-cycle$({\bf u}_{l-1}, {\bf s}_{l-1},l-1)$ (Computation of the coarse grid correction)
\item ${\bf u}_l = {\bf u}_l + {\bf P}_{l,l-1}({\bf u}_{l-1}-{\bf \tilde{u}}_{l-1})$ (Correction via Prolongation)
\end{itemize}
\item end if
\end{itemize}

The restriction ${\bf R}_{l-1,l}$ is an agglomeration that weighs components by the volume of their cells and divides by the total volume. As for the prolongation ${\bf P}_{l,l-1}$, it uses a bilinear weighting \cite{jameso:86}.

On the finest level, the smoother is applied to the equation \eqref{nonlin-eq} resp. \eqref{nonlin-eq-steady}. On sublevels, it is instead used to solve
\begin{equation}\label{coarse-nonlin-eq}
{\bf F}:={\bf s}_l-{\bf F}_l({\bf u}_l) = {\bf 0}
\end{equation}
with
\[{\bf s}_l={\bf F}_l({\bf R}_{l,l+1} {\bf u}_{l+1}) + {\bf R}_{l,l+1}{\bf r}_{l+1}.\]

\section{Preconditioned smoothers}

All smoothers we use have a pseudo time iteration as a basis. These are iterative methods for the nonlinear equation ${\bf F}({\bf u})={\bf 0}$ that are obtained by applying a time integration method to the initial value problem
\[{\bf u}_{t^*} =- {\bf F}({\bf u}), \quad {\bf u}(0) = {\bf u}^0.\]
For convenience, we have dropped the subscript $l$ that denotes the multigrid level.

\subsection{Preconditioned additive Runge-Kutta methods}

We start with splitting ${\bf F}({\bf u})$ in a convective and diffusive part
\begin{equation}\label{split-rhs}
{\bf F}({\bf u})={\bf f}^c({\bf u}) + {\bf f}^v({\bf u}).
\end{equation}
Hereby, ${\bf f}^c$ contains the physical convective fluxes, as well as the discretized time derivative and the multigrid source terms, whereas ${\bf f}^v$ contains both the artificial dissipation and the discretized second order terms of Navier-Stokes. 

An additive explicit Runge-Kutta (AERK) method is then implemented in the following form:
\begin{align}
{\bf u}^{(0)} & = {\bf u}\\
{\bf u}^{(i)} & = {\bf u} - \alpha_{i} \Delta t^*({\bf f}^{c,(i-1)} + {\bf f}^{v,(i-1)}), \quad i=1,...,s \\
{\bf u}^{n+1} & = {\bf u}^{(s)},
\end{align}
where
\begin{align}
{\bf f}^{c,(i)} & = {\bf f}^{c}({\bf u}^{(i)}), \quad i=0,...,s-1 \label{stage derivatives 1}\\
{\bf f}^{v,(0)} & = {\bf f}^{v}({\bf u}^{(0)}), \label{stage derivatives 2}\\
{\bf f}^{v,(i)} & = \beta_{j+1}{\bf f}^{v}({\bf u}^{(i)}) + (1-\beta_{i+1}){\bf f}^{v,(i-1)}, \quad i=1,...,s-1. \label{stage derivatives 3}
\end{align}
The second to last line implies that $\beta_1=1$. Here, $\Delta t^*$ is a local pseudo time step, meaning that it depends on the specific cell and the multigrid level. It is obtained by choosing $c^*$, a CFL number in pseudo time, and then computing $\Delta t^*$ based on the local mesh width $\Delta x_{k_l}$:
\[\Delta t^* = c^* \Delta x_{k_l}\]
This implies larger time steps on coarser cells, in particular on coarser grids. 

\begin{table}[h]
\centering
\begin{tabular}{cc|ccccc}
& $i$ & 1 & 2 & 3 & 4 & 5 \\ \hline
ARK3J & $\alpha_i$ & 0.1481 & 2/5 & 1 & - & - \\
 & $\beta_i$  & 1 & 1/2 & 1/2 & - & - \\\hline
ARK5J & $\alpha_i$ & 1/4 & 1/6 & 3/8 & 1/2 & 1 \\
$\beta_i$ && 1 & 0 & 0.56 & 0 & 0.44 \\\hline\hline
ARK51 & $\alpha_i$ & 0.0533 & 0.1263 & 0.2375 & 0.4414& 1. \\ \hline
ARK52 & $\alpha_i$ & 0.0695 & 0.1602 & 0.2898 & 0.5060 & 1. 
\end{tabular}
\caption{\label{RK smoothers}Coefficients of explicit and additive Runge-Kutta smoothers, 3- and 5-stage method.}
\end{table}
As for the coefficients, 3-, 4- and 5-stage schemes have been designed to have good smoothing properties in a multigrid method for convection dominated model equations. The first three schemes have been designed by Jameson using linear advection with a fourth order diffusion term. We denote these by ARKsJ with $s$ the number of stages. See \cite{jameso:85} for ARK4J and ARK5J and \cite{jameso:15} for ARK3J. The 5-stage schemes ARK51 and ARK52 are from \cite{letapo:89}. The scheme ARK52 is employed in \cite{swturo:07}. Coefficients for the 3- and 5-stage schemes can be found in table \ref{RK smoothers}. All of these schemes are first order, except for the last one, which has order two and is therefore denoted as ARK52. In the original publication ARK51 and ARK52 are not additive. When using these within an additive method, we use the $\beta$ coefficients from ARK5J. For current research into improving these coefficients we refer to \cite{birken:12,bedods:17}.

Setting $\beta_j=1$ for all $j$ gives an unsplit low storage explicit Runge-Kutta method that does not treat convection and diffusion differently. We refer to these schemes as ERK methods, e.g. ERK3J or ERK51. 

To precondition this scheme, a preconditioner ${\bf P}^{-1} \in \mathbb{R}^{m\times m}$ is applied to the equation system \eqref{nonlin-eq} or \eqref{coarse-nonlin-eq} by multiplying them with it, resulting in an equation
\[{\bf P}^{-1}{\bf F}({\bf u})={\bf 0}.\]
In a pseudo-time iteration for the new equation, all function evaluations have to be adjusted. In the above algorithm, this is realized by replacing the term $\alpha_i \Delta t^*({\bf f}^{c,(i-1)} + {\bf f}^{v,(i-1)})$ with $\alpha_i \Delta t^*{\bf P}^{-1}({\bf f}^{c,(i-1)} + {\bf f}^{v,(i-1)})$. 

A good preconditioner should approximate the Jacobian $\frac{\partial {\bf F}}{\partial {\bf u}}$ of ${\bf F}$ well, while simultaneously being easy to apply. 






\subsection{Additive W-methods}

An alternative way of deriving a preconditioned explicit method has been presented by Langer in \cite{langer:11}. He calls these methods preconditioned implicit smoothers and derives them from specific singly diagonally implicit RK (SDIRK) methods. SDIRK methods consist of a nonlinear system at each step, which he solves with one Newton step each and then simplifies by always using the Jacobian from the first stage. This is known as a special choice of Rosenbrock method in the literature on differential equations \cite[p. 102]{haiwan:04}. To arrive at a preconditioned method, Langer then replaces the system matrix with an approximation, for example originating from a preconditioner as known from linear algebra. In fact, this type of method is called a W-method in the IVP community \cite[p. 114]{haiwan:04}. 

\begin{table}[h]
\centering
\begin{tabular}{c|cccc}
 & $\eta$ & 0 & 0 & 0\\
 & $\alpha_{1}$ & $\eta$ & 0 & 0\\
 & 0 & $\ddots$ & $\ddots$ & 0 \\ 
 & 0 & 0 & $\alpha_{s-1}$ & $\eta$ \\ \hline
 & 0 & $\ldots$ & 0 & $\alpha_s$
\end{tabular}
\caption{\label{ASDIRK-coeff-c}Butcher arrays for additive SDIRK method: Convective terms.}
\end{table}

\begin{table}[h]
\centering
\begin{tabular}{c|ccccc}
& $\eta$ &  & $\ldots$ & & 0 \\
& $\alpha_{1}$ & $\eta$ &  &  & 0\\
& $\alpha_2(1-\beta_1)$ & $\alpha_2\beta_2$ & $\ddots$ &  & 0\\ 
& 0 & $\ddots$ & $\ddots$ & $\eta$ & 0\\ 
& 0 & $\ldots$ & $\alpha_{s-1}(1-\beta_{s-1})$ & $\alpha_{s-1}\beta_s$ & $\eta$\\ \hline
& 0 & $\ldots$ & 0 & $\alpha_s(1-\beta_{s-1})$ & $\alpha_s\beta_s$ 
\end{tabular}
\caption{\label{ASDIRK-coeff-v}Butcher arrays for additive SDIRK method: Diffusive terms.}
\end{table}

We now extend the framework from \cite{langer:11} to additive Runge-Kutta methods. For clarity we repeat the derivation, but start from the split equation
\begin{equation}\label{split-eq}
{\bf u}_{t^*} + {\bf f}^c({\bf u}) + {\bf f}^v({\bf u}) = {\bf 0}
\end{equation}
as described in \eqref{split-rhs}. To this equation, we apply an additive SDIRK method with coefficients given in tables \ref{ASDIRK-coeff-c} and \ref{ASDIRK-coeff-v}:

\begin{align}
{\bf k}_i & =-{\bf F}({\bf u}^n  + \Delta t^*(\sum_{j=1}^{i-1}(a^c_{ij}{\bf k}^c_j + a^v_{ij}{\bf k}^v_j)+ \eta{\bf k}_i)), \quad i=1,...,s, \label{additive-SDIRK} \\
{\bf u}^{n+1} & = {\bf u}^n + \Delta t^* (\alpha_s{\bf k}^c_s + \alpha_s(1-\beta_{s-1}){\bf k}^v_{s-1} + \alpha_s\beta_s{\bf k}^v_s). \label{newvalue}
\end{align}
Hereby, the vectors ${\bf k}$ are called stage derivatives and we have ${\bf k} = {\bf k}^c + {\bf k}^v$ according to the splitting \eqref{split-rhs}. Thus, we have to solve $s$ nonlinear equation systems for the stage derivatives ${\bf k}_i$. 

To obtain an additive Rosenbrock method, these are solved approximately using one Newton step each with initial guess zero, changing the stage values to

\begin{align}\label{additive Rosenbrock}
{\bf k}_i & =-({\bf I}+\eta\Delta t^*{\bf J}_i)^{-1}{\bf F}({\bf u}^n  + \Delta t^*(\sum_{j=1}^{i-1}(a^c_{ij}{\bf k}^c_j + a^v_{ij}{\bf k}^v_j))), \quad i=1,...,s, 
\end{align}
where ${\bf J}_i = \frac{\partial {\bf F}_i({\bf 0})}{\partial{\bf k}}$, with ${\bf F}_i({\bf k}) := {\bf F}({\bf u}^n  + \Delta t^*(\sum_{j=1}^{i-1}(a^c_{ij}{\bf k}^c_j + a^v_{ij}{\bf k}^v_j)+ \eta{\bf k}))$. Thus, we now have to solve a linear system at each stage. This type of scheme is employed in Swanson et. al. \cite{swturo:07}. They refer to the factor $\eta$ as $\epsilon$ and provide a discrete Fourier analysis of this factor. 

As a final step, we approximate the system matrices ${\bf I}+\eta\Delta t^*{\bf J}_i$ by a matrix ${\bf W}$. This gives us a new class of schemes, which we call additive W (AW) methods, with stage derivatives given by:

\begin{align}\label{additive W}
{\bf k}_i & =-{\bf W}^{-1}{\bf F}({\bf u}^n  + \Delta t^*(\sum_{j=1}^{i-1}(a^c_{ij}{\bf k}^c_j + a^v_{ij}{\bf k}^v_j))), \quad i=1,...,s,
\end{align}

In both additive Rosenbrock and additive W methods, equation \eqref{newvalue} remains unchanged. 

Finally, after some algebraic manipulations, this method can be rewritten in the same form as the low storage preconditioned ARK methods presented earlier:
\begin{align*}
{\bf u}^{(0)} & = {\bf u}^n\\
{\bf u}^{(i)} & = {\bf u}^n - \alpha_{i} \Delta t^*{\bf W}^{-1}({\bf f}^{c,(i-1)} + {\bf f}^{v,(i-1)}), \quad i=1,...,s \\
{\bf u}^{n+1} & = {\bf u}^{(s)},
\end{align*}
where
\begin{align*}
{\bf f}^{c,(i)} & = {\bf f}^{c}({\bf u}^{(i)}), \quad i=0,...,s-1 \\
{\bf f}^{v,(0)} & = {\bf f}^{v}({\bf u}^{(0)}),\\
{\bf f}^{v,(i)} & = \beta_{j+1}{\bf f}^{v}({\bf u}^{(i)}) + (1-\beta_{i+1}){\bf f}^{v,(i-1)}, \quad i=1,...,s-1.
\end{align*}

As for the explicit methods, one recovers an unsplit scheme for $\beta_i=1$ for all $i$ and we refer to these methods as SDIRK, Rosenbrock and W methods.

\subsection{Comparison}
\label{comparison}
To get a better understanding for the different methods, it is illustrative to consider the linear case. Then, these methods are iterative schemes to solve a linear equation system $({\bf A} + {\bf B}){\bf x} = {\bf b}$ and can be written as
\[{\bf x}^{k+1} = {\bf Mx}^k + {\bf Nb}.\]
The matrix ${\bf M}$ is the iteration matrix and for pseudo time iterations, it is given as the stability function $S$ of the time integration method. These are a polynomial $P_s$ of degree $s$ in $\Delta t^*({\bf A} + {\bf B})$ for an $s$ stage ERK method and a bivariate polynomial $P_s$ of degree $s$ in $\Delta t^*{\bf A}$ and $\Delta t^*{\bf B}$ for an $s$ stage AERK method. When preconditioning is added, this results in $P_s(\Delta t^*{\bf P}^{-1}({\bf A}+{\bf B}))$ and $P_s(\Delta t^*{\bf P}^{-1}{\bf A},\Delta t^*{\bf P}^{-1}{\bf B})$, respectively. 

For the implicit schemes, we obtain a rational function of the form $Q_s({\bf I}+\eta\Delta t^*({\bf A} + {\bf B}))^{-1}P_s(\Delta t^*({\bf A} + {\bf B}))$ for an $s$ stage SDIRK or Rosenbrock method. Here, $Q_s$ is a second polynomial of degree $s$. Finally, for an $s$-stage additive W method, we obtain a function of the form $Q_s({\bf W})^{-1}P_s(\Delta t^*{\bf A},\Delta t^*{\bf B})$. Due to the specific contruction, the inverse can simply be moved from the left into the argument which gives $P_s(\Delta t^*{\bf W}^{-1}{\bf A},\Delta t^*{\bf W}^{-1}{\bf B})$. Note that this is the same as for the preconditioned AERK method, except for the preconditioner. 

The additive W method and the preconditioned ARK method have three main differences. First of all, there is the role of ${\bf P}$ in the AERK method versus the matrix ${\bf W}$. In the W method ${\bf W}\approx ({\bf I}+\eta\Delta t^*{\bf J}_i)$, whereas in the AERK scheme, ${\bf P} \approx {\bf J}_i$. Second, the timestep in the one case is that of an explicit ARK method, whereas in the other, that of an implicit method. The latter in its SDIRK or Rosenbrock form is A-stable. However, approximating the Jacobian can cause the stability region to become finite. Finally, the latter method has an additional parameter $\eta$ that needs to be chosen. However, the large stability region makes the choice of $\Delta t^*$ easy for the additive W method (very large), whereas it has to be a small value for the preconditioned ARK scheme.

\subsection{SGS Preconditioner}

The basis of our method is the preconditioner suggested by Swanson et al. in \cite{swturo:07} and improved by Jameson in \cite{jameso:17}. In effect, this is a choice of a ${\bf W}$ matrix in the framework just presented. We now repeat the derivation of their preconditioner in our notation to obtain an improved version. 

The first step is to approximate the Jacobian by using a different first order linearized discretization. It is based on a splitting ${\bf A} = {\bf A}^+ + {\bf A}^-$ of the flux Jacobian. This is evaluated in the average of the values on both sides of the interface, thereby deviating from \cite{swturo:07}. The split Jacobians correspond to positive and negative eigenvalues:
\[{\bf A}^+ = \frac12 ({\bf A} + |{\bf A}|), \quad  {\bf A}^- = \frac12 ({\bf A} - |{\bf A}|).\]
Alternatively, these can be written in terms of the matrix of right eigenvectors ${\bf R}$ as
\[{\bf A}^+ = {\bf R}|\Lambda^+|{\bf R}^{-1}, \quad  {\bf A}^- = {\bf R}|\Lambda^-|{\bf R}^{-1},\]
where $\Lambda^{\pm}$ are diagonal matrices containing the positive and negative eigenvalues, respectively. 

As noted in \cite{jameso:17}, it is now crucial to use a cutoff function for the eigenvalues beforehand, to bound them away from zero. We use a parabolic function which kicks in when the modulus of the eigenvalue $\lambda$ is smaller or equal to a fraction $ad$ of the speed of sound $a$ with free parameter $d \in [0,1]$:
\begin{equation}\label{cutoff}
|\lambda| = \frac{1}{2}\left(ad + \frac{|\lambda|^2}{ad} \right), \quad |\lambda| \leq ad.
\end{equation}

With this, an upwind discretization is given in cell $i$ by
\begin{equation}\label{steadylinearized}
{\bf u}_{i_t} = \frac{1}{\Omega_i}\sum_{e_{ij}\in N(i)}|e_{ij}|({\bf A}_{{\bf n}_{ij}}^+{\bf u}_i + {\bf A}_{{\bf n}_{ij}}^-{\bf u}_j).
\end{equation}
Here, $e_{ij}$ is the edge between cells $i$ and $j$, $N(i)$ is the set of cells neighboring $i$ and ${\bf n}_{ij}$ the unit normal vector from $i$ to $j$. 

For the unsteady equation \eqref{nonlin-eq}, we obtain instead
\begin{equation}\label{form1}
{\bf u}_{i_t}=\frac{3}{2\Delta t^*} {\bf I} + \frac{1}{\Omega_i}\sum_{e_{ij}\in N(i)}|e_{ij}|({\bf A}^+{\bf u}_i + {\bf A}^-{\bf u}_j).
\end{equation}

The corresponding approximation of the Jacobian is then used to construct a preconditioner. Specifically, we consider the block SGS preconditioner
\begin{equation}\label{sgs-precon}
{\bf P}^{-1} = ({\bf D} + {\bf L}) {\bf D}^{-1} ({\bf D} +{\bf U}),
\end{equation}
where ${\bf L}$, ${\bf D}$ and ${\bf U}$ are block matrices with $4\times 4$ blocks. This preconditioner would look different when several SGS steps would be performed. However, we did not find this to be beneficial. 

We now have two cases. In the AERK framework, ${\bf L} + {\bf D} + {\bf U} = {\bf J}$ and we arrive at 
\begin{eqnarray}\label{d-l-u-matrices-J}
{\bf L}_{ij} = -\frac{1}{\Omega_{i}}
  (\Delta y {\bf A}_{i-1,j}^+ + \Delta x {\bf B}_{i,j-1}^+ ),\\
{\bf U}_{ij} = \frac{1}{\Omega_{i}}
  (\Delta y {\bf A}_{i-1,j}^- + \Delta x {\bf B}_{i,j-1}^- ),\\
{\bf D}_{ii} = \frac{1}{\Omega_i} [ \Delta y ({\bf A}_{ii}^+ -{\bf A}_{ii}^-) + \Delta x({\bf B}_{ii}^+ -{\bf B}_{ii}^-)],
\end{eqnarray}
respectively
\begin{equation}\label{d-unsteady-J}
{\bf D}_{ii} = \frac{3}{2\Delta t} {\bf I} +\frac{1}{\Omega_i} [ \Delta y ({\bf A}_{ii}^+ -{\bf A}_{ii}^-) + \Delta x({\bf B}_{ii}^+ -{\bf B}_{ii}^-)],
\end{equation}
in the unsteady case. 

In the additive W framework, ${\bf L} + {\bf D} + {\bf U} = {\bf I} + \eta\Delta t^*{\bf J}$ and we obtain
\begin{eqnarray}\label{d-l-u-matrices-W}
{\bf L}_{ij} = -\frac{\eta\Delta t_i^*}{\Omega_{i}}
  (\Delta y {\bf A}_{i-1,j}^+ + \Delta x {\bf B}_{i,j-1}^+ ),\\
{\bf U}_{ij} = \frac{\eta\Delta t_i^*}{\Omega_{i}}
  (\Delta y {\bf A}_{i-1,j}^- + \Delta x {\bf B}_{i,j-1}^- ),\\
{\bf D}_{ii} = {\bf I} + \frac{\eta\Delta t_i^*}{\Omega_i} [ \Delta y ({\bf A}_{ii}^+ -{\bf A}_{ii}^-) + \Delta x({\bf B}_{ii}^+ -{\bf B}_{ii}^-)].
\end{eqnarray}
or in the unsteady case
\begin{equation}\label{d-unsteady-W}
{\bf D}_{ii} = {\bf I} + \frac{3\eta\Delta t^*}{2\Delta t} {\bf I} + \frac{\eta\Delta t_i^*}{\Omega_i} [ \Delta y ({\bf A}_{ii}^+ -{\bf A}_{ii}^-) + \Delta x({\bf B}_{ii}^+ -{\bf B}_{ii}^-)].
\end{equation}

Applying this preconditioner requires solving small $4\times 4$ systems coming from the diagonal. We use Gaussian elimination for this. A fast implementation is obtained by transforming first to a certain set of symmetrizing variables, see \cite{swturo:07}.

\section{Discrete Fourier Analysis}

We now perform a discrete Fourier analysis of the preconditioned ARK method for the two dimensional Euler equations using the JST scheme. For a description of this technique, also called local Fourier analysis (LFA) in the multigrid community, we refer to \cite{troosc:01,gustaf:08}. The rationale for this is that the core convergence problems for multigrid methods for viscous flow problems on high aspect ratio grids are the convective terms and the high aspect ratio grids. The viscous terms are of comparatively minor importance. Here, we do not take into account the coarse grid correction. Thus, our aim is to obtain amplification- and smoothing factors for the smoother. The latter is known to be representative of 2-grid convergence rates and is given by
\begin{equation}\label{smoothfactor}
\max_{\lambda_{HF}}|S(\lambda)|,
\end{equation}
where $\lambda_{HF}$ denote the high frequency eigenvalues. Since eigenfunctions of first order hyperbolic differential operators involve $e^{i\phi x}$, these are in $[-\pi,-\pi/2]$ and $[\pi/2,\pi]$. 

We now consider a linearized version of the underlying equation with periodic boundary conditions on the domain $\Omega = [0,1]^2$:
\begin{equation}\label{linearized-euler}
\frac{d}{dt}{\bf u} + ({\bf A}{\bf u})_x + ({\bf B}{\bf u})_y = {\bf 0}
\end{equation}
with ${\bf A} = \frac{\partial {\bf f}_1}{\partial {\bf u}}$ and ${\bf B} = \frac{\partial {\bf f}_2}{\partial {\bf u}}$ being the Jacobians of the Euler fluxes in a fixed point $\hat{{\bf u}}$, to be set later. 

\subsection{JST scheme}
We discretize \eqref{linearized-euler} on a cartesian mesh with mesh width $\Delta x$ in $x$-direction and $\Delta y = \text{AR}\Delta x$ in $y$-direction (AR=aspect ratio), resulting in an $n_x\times n_y$ mesh. A cell centered finite volume method with the JST flux is employed. We denote the shift operators in $x$ and $y$ direction by $E_x$ and $E_y$. Cells are indexed the canonical doubly lexicographical way for a cartesian mesh. In cell $ij$  we write the discretization as
\[({\bf H}{\bf u})_{ij} = (({\bf H}_c+{\bf H}_v){\bf u})_{ij}\]
with
\[{\bf H}_c = \frac{1}{2\Delta x\Delta y}({\bf A}(E_x^{+1}-E_x^{-1})\Delta y + {\bf B}(E_y^{+1}-E_y^{-1})\Delta x),\]
respectively in the unsteady case, 
\[{\bf H}_c = \frac{3}{2 \Delta t}{\bf I} + \frac{1}{2\Delta x \Delta y}({\bf A}(E_x^{+1}-E_x^{-1})\Delta y + {\bf B}(E_y^{+1}-E_y^{-1})\Delta x).\]

For ${\bf H}_v$, the starting point is that the pressure in conservative variables is 
\[p = (\gamma -1)\left( \rho E-\rho\frac{(\rho v_1)^2+(\rho v_2)^2}{2\rho^2} \right ).\]
In the fraction, all potential shift operators cancel out. Thus, for the second order differences in both directions, 
\[p_{j+1}-2p_j+p_{j-1} = (\gamma -1)[(E^+-2+E^-)\rho E_j - |{\bf v}|^2/2(E^+-2+E^-)\rho_j].\]
For the fourth order difference, there's a corresponding identity. Furthermore, applying the second or fourth order difference to $\rho H_j = \rho E_j + p_j$ results in 
\[\rho H_{j+1}-2\rho H_j+\rho H_{j-1} = \gamma (E^+-2+E^-)\rho E_j - (\gamma -1)|{\bf v}|^2/2(E^+-2+E^-)\rho_j.\] 
This gives
\begin{eqnarray*}
{\bf H}_v = \frac{1}{\Delta x \Delta y}{\bf M}[\epsilon^{(2)}((-E_x^{+1}+2-E_x^{-1})\Delta y + (-E_y^{+1}+2-E_y^{-1})\Delta x) \\
 +\epsilon^{(4)}(E_x^{+2} - 4E_x^{+1} + 6 - 4E_x^{-1}+E_x^{-2})\Delta y\\
 +\epsilon^{(4)}(E_y^{+2} - 4E_y^{+1} + 6 - 4E_y^{-1}+E_y^{-2})\Delta x]
\end{eqnarray*}
with
\[{\bf M} = \left( \begin{array}{cccc} 1 & 0 & 0 & 0\\
0 & 1 & 0 & 0 \\
0 & 0 & 1 & 0 \\
-(\gamma -1)|{\bf v}|^2/2 & 0 & 0 & \gamma 
\end{array}\right ).\]

For the coefficient functions $\epsilon^{(2)}$ and $\epsilon^{(4)}$ (see \eqref{e2} and \eqref{e4}), we first look at the shock sensor $s_{j+1/2}$. Here, we use the version for the Euler equations based on pressure. Straightforward calculations give
\[p_{j+1}+2p_j+p_{j-1} = (\gamma -1)[(E^{+1}+2+E^{-1})\rho E_j-2|{\bf v}|^2(E^++2+E^-)\rho_j].\]
Thus
\[s_j = \frac{(\gamma -1)[(E^{+1}-2+E^{-1})\rho E_j-1/2|{\bf v}|^2(E^+-2+E^-)\rho_j]}{(\gamma -1)[(E^{+1}+2+E^{-1})\rho E_j-2|{\bf v}|^2(E^++2+E^-)\rho_j] + 0.001}.\]
For simplicity, we now assume that $\max(s_j,s_{j+1}) = s_j$. Thus,
\[s_{j+1/2} = \min(0.25,s_j).\]

For the spectral radius we note that in the speed of sound $a_j=\sqrt{\gamma p_j/\rho_j}$, possible shift operators cancel out as well, implying that is constant over the mesh. This gives
\[r_i = |v_1|+a,\]
\[r_j = |v_2|+a.\]
Regarding the maxima, we have $r_j=r_{j+1} =: r$ and correspondingly for the $y$ direction with $r_i$. Thus, 
\[\epsilon^{(2)} = rs_{j+1/2}\]
and
\[\epsilon^{(4)} = \max(0,r/32-2\epsilon^{(2)}).\]

\subsection{Preconditioner}

With regards to the SGS preconditioner ${\bf P}^{-1} = ({\bf D} + {\bf L}) {\bf D}^{-1} ({\bf D} +{\bf U})$, the different discretization based on the flux splitting \eqref{steadylinearized} with cutoff function \eqref{cutoff} gives (see \eqref{d-l-u-matrices-J}-\eqref{d-unsteady-J})
\[ {\bf L} = -\frac{1}{\Omega} [\Delta y{\bf A}^+ E_x^{-1} + \Delta x{\bf B}^+E_y^{-1}], \]
\[ {\bf U} = \frac{1}{\Omega} [\Delta y{\bf A}^- E_x^{+1} + \Delta x{\bf B}^-E_y^{+1}]. \]
We now get two different operators for the diagonal part for the steady and for the unsteady case. We have
\[{\bf D}^s = {\bf I} + \frac{1}{\Omega} [\Delta y({\bf A}^+-{\bf A}^-)+(\Delta x ({\bf B}^+-{\bf B}^-)],\]
for the steady case, whereas for the unsteady case there is
\[{\bf D}^u = \frac{3}{2\Delta t}{\bf I} + {\bf D}^s.\]
With these, the preconditioner \eqref{sgs-precon} is formed. For the W methods, these matrices need to be adjusted slightly, compare \eqref{d-l-u-matrices-W}-\eqref{d-unsteady-W}. 

We now make one simplification in the analysis and that is that we assume the matrices to be evaluated with the value of the respective cell and not the average as in the actual method. 

As an example, the application of the 3-stage ARK scheme results in the following operator, where we write $\bar{{\bf H}}_c:={\bf P}^{-1}{\bf H}_c$, $\bar{{\bf H}}_c:={\bf P}^{-1}{\bf H}_v$ and $\bar{\alpha}_i = \Delta t^*\alpha_i$:

\begin{eqnarray*}
{\bf G} = {\bf I} -\bar{\alpha}_3((\bar{{\bf H}}_c + \beta_3\bar{{\bf H}}_v)({\bf I} - \bar{\alpha}_2((\bar{{\bf H}}_c + \beta_2\bar{{\bf H}}_v)({\bf I}-\bar{\alpha}_1(\bar{{\bf H}}_c + \bar{{\bf H}}_v)) + (1-\beta_2)\bar{{\bf H}}_v)) \nonumber \\
+ (1-\beta_3)(\beta_2\bar{{\bf H}}_v(I-\bar\alpha_1(\bar{{\bf H}}_c + \bar{{\bf H}}_v)) + (1-\beta_2)\bar{{\bf H}}_v)).
\end{eqnarray*}
For other smoothers, we have to use other appropriate stability functions, as discussed in section~\ref{comparison}.

\subsection{Amplification and Smoothing factors}
We are now interested in the amplification factor of the corresponding method for different values of $\Delta x$ and $\Delta y$. Working with ${\bf G}$ directly would require assembling a large matrix in $\mathbb{R}^{4n_x\times 4n_y}$. Instead, we perform a discrete Fourier transform. In Fourier space, the transformed operator block diagonalizes, allowing to work with the much smaller matrix ${\bf \hat{G}}\in\mathbb{R}^{4\times4}$. Thus, we replace ${\bf u}_{ij}$ by its discrete Fourier series
\[{\bf u}_{ij} = \sum_{k_x=-n_x/2+1}^{n_x/2}\sum_{k_y=-n_y/2+1}^{n_y/2}{\bf \hat{u}}_{k_x,k_y}e^{2\pi i(k_xx_i + k_yy_j)}\]
and analyze 
\[{\bf \hat{u}}_{k_x,k_y}^{k+1} = {\bf \hat{G}}_{k_x,k_y}{\bf \hat{u}}_{k_x,k_y}^k.\]
When applying a shift operator to one of the exponentials, we obtain
\[E_xe^{2\pi i(k_xx_i + k_yy_j)} =e^{2\pi i(k_x(x_i+1/n_x) + k_yy_j)} = e^{2\pi ik_x/n_x}e^{2\pi i(k_xx_i + k_yy_j)}\]
and similar for $E_y$. Defining the phase angles
\[\Theta_x = 2\pi k_x/n_x, \quad \Theta_y = 2\pi k_y/n_y,\]
the Fourier transformed shift operators are
\[\hat{E}_x = e^{i\Theta_x}, \quad \hat{E}_y = e^{i\Theta_y}\]
and can replace the dependence on the wave numbers with a dependence on phase angles. 

To compute the spectral radius of ${\bf G}$, we now just need to look at the maximum of the spectral radius of ${\bf \hat{G}}_{\Theta_x,\Theta_y} = {\bf \hat{G}}_{k_x,k_y}$ over all phase angles $\Theta_x$ and $\Theta_y$ between $-\pi$ and $\pi$. Furthermore, this allows to compute the smoothing factor \eqref{smoothfactor} as well, by instead taking the maximum over all wave numbers between $-\pi$ and $-\pi/2$, as well as $\pi/2$ and $\pi$.

\subsection{Results}

We evaluate the matrices in the points 
\[\hat{{\bf u}}_1 =(1,\sqrt{2}/2,\sqrt{2}/2,3.290)^T \quad\text{(Mach 0.8, $\alpha = 45^{\circ}$)}\]
and
\[\hat{{\bf u}}_3 =(1,1,0,3.290)^T \quad\text{(Mach 0.8, $\alpha = 0^{\circ}$)}.\]
We use a $8 \times (8\cdot AR)$ grid with different aspect ratios (AR), namely AR=1, AR=100 and AR=10000. To determine the physical time step, a CFL number $c$ of 200 is chosen. All results were obtained using a python script, which can be accessed at http://www.maths.lu.se/philipp-birken/rksgs\_fourier.zip.

\subsubsection{The explicit schemes}

\begin{table}
\centering
\begin{tabular}{c|ccc|ccc}
& \multicolumn{3}{c|}{ERK3} & \multicolumn{3}{c}{ARK3J} \\
AR              & 1 & 100 & 10000 & 1 & 100 & 10000\\ \hline
$\rho({\bf M})$ & 0.9933 & 0.9935 & 0.9935 & 0.9933 & 0.9935 & 0.9935 \\
Sm. fct.        & 0.5158 & 0.9935 & 0.9935 & 0.4634 & 0.9935 & 0.9935 \\ \hline
\end{tabular}
\caption{\label{ERK-ARK}Amplification and smoothing factors of ERK3 and ARK3J, $8 \times (8AR)$ grid, $c$=0.9; $M=0.8$, $\alpha = 0^{\circ}$}
\end{table}

Results for explicit schemes for different test cases are shown in table \ref{ERK-ARK}. As can be seen, these methods have terrible convergence rates, but are good smoothers for equidistant meshes. For non-equidistant meshes, this is not the case, which demonstrates the poor performance of these methods for viscous flow problems. 

\subsubsection{Preconditioned ARK}

\begin{figure}[h]
\centering
\includegraphics[width=0.49\textwidth]{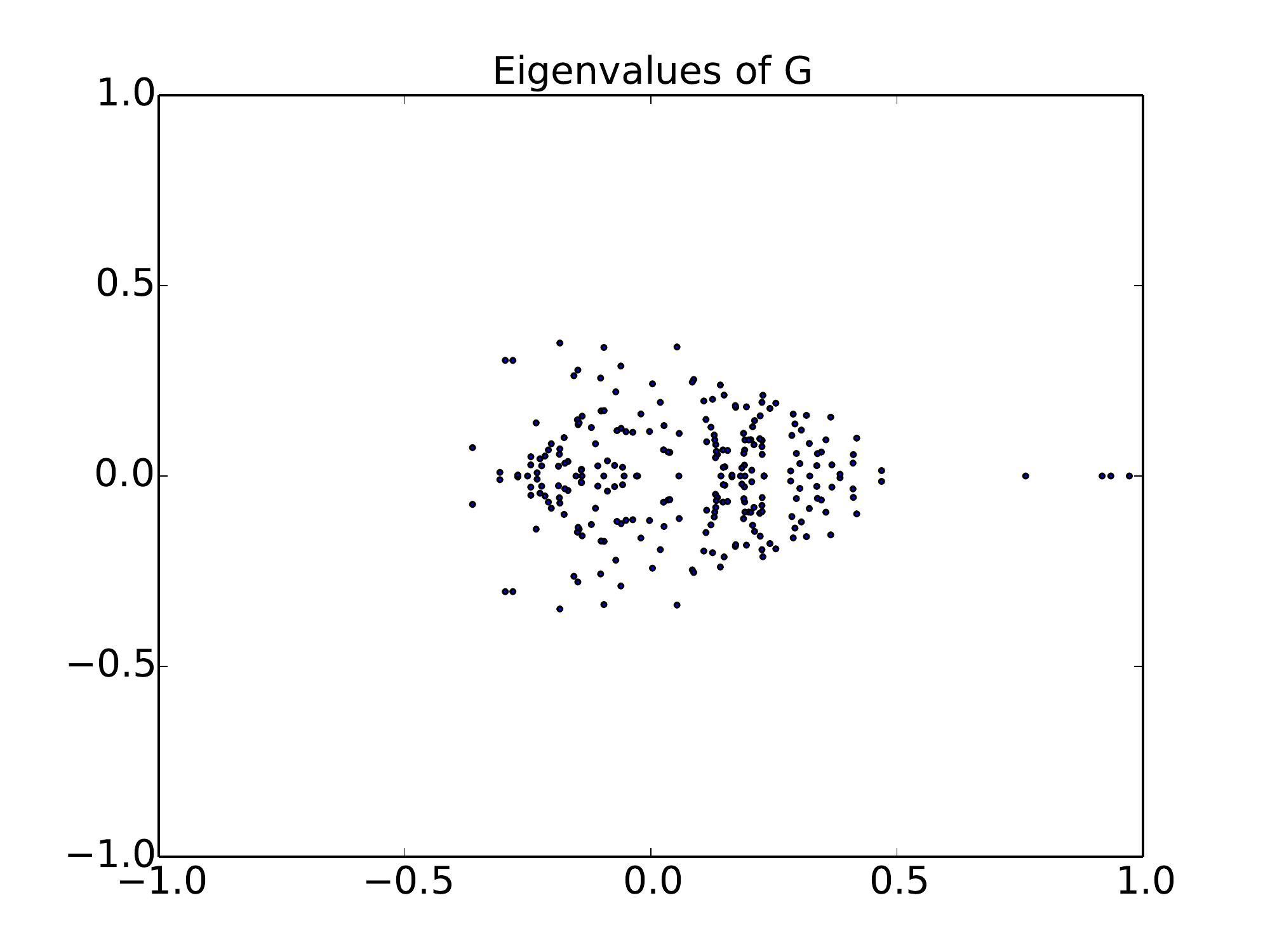}
\includegraphics[width=0.49\textwidth]{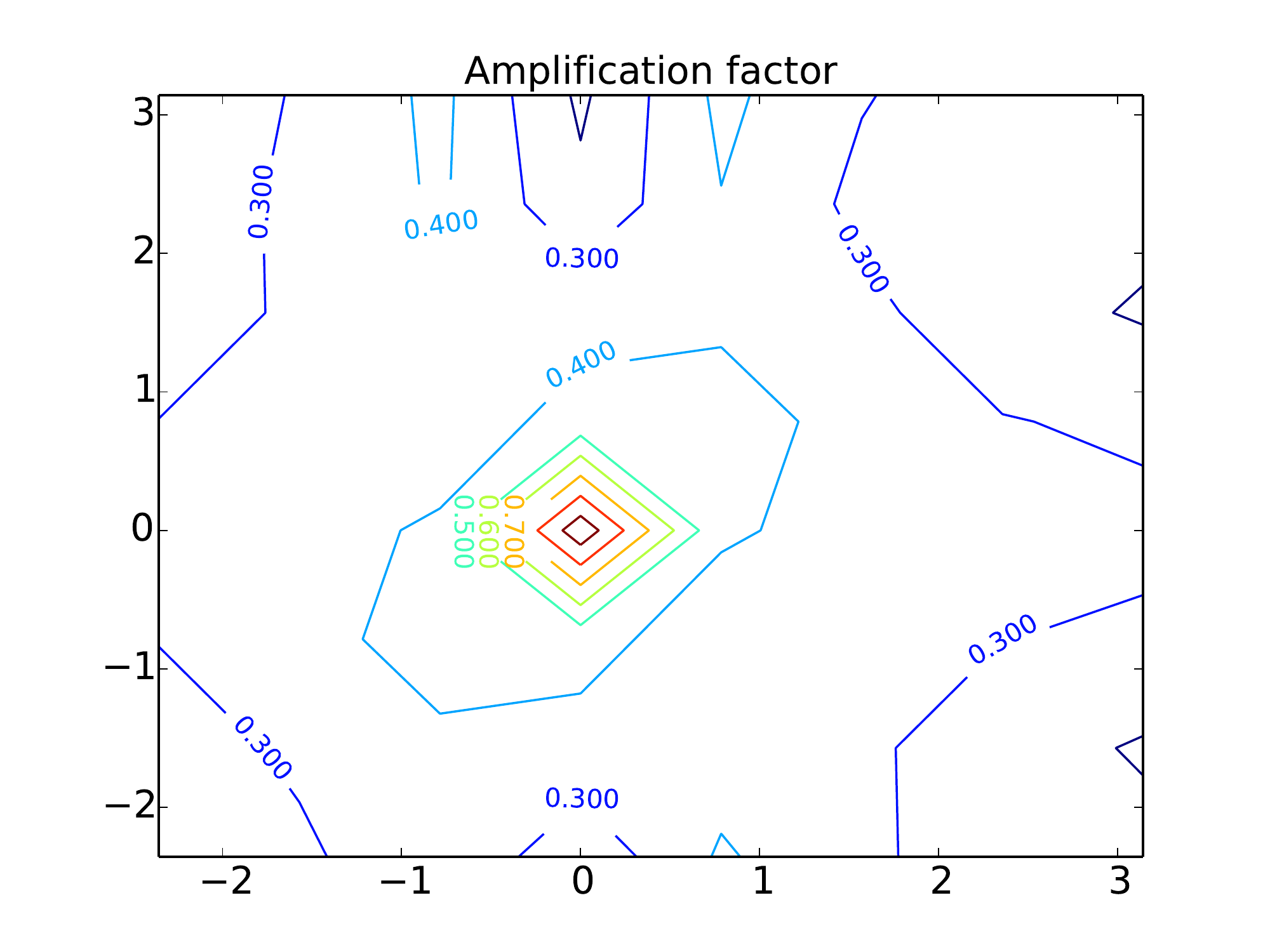}
\includegraphics[width=0.49\textwidth]{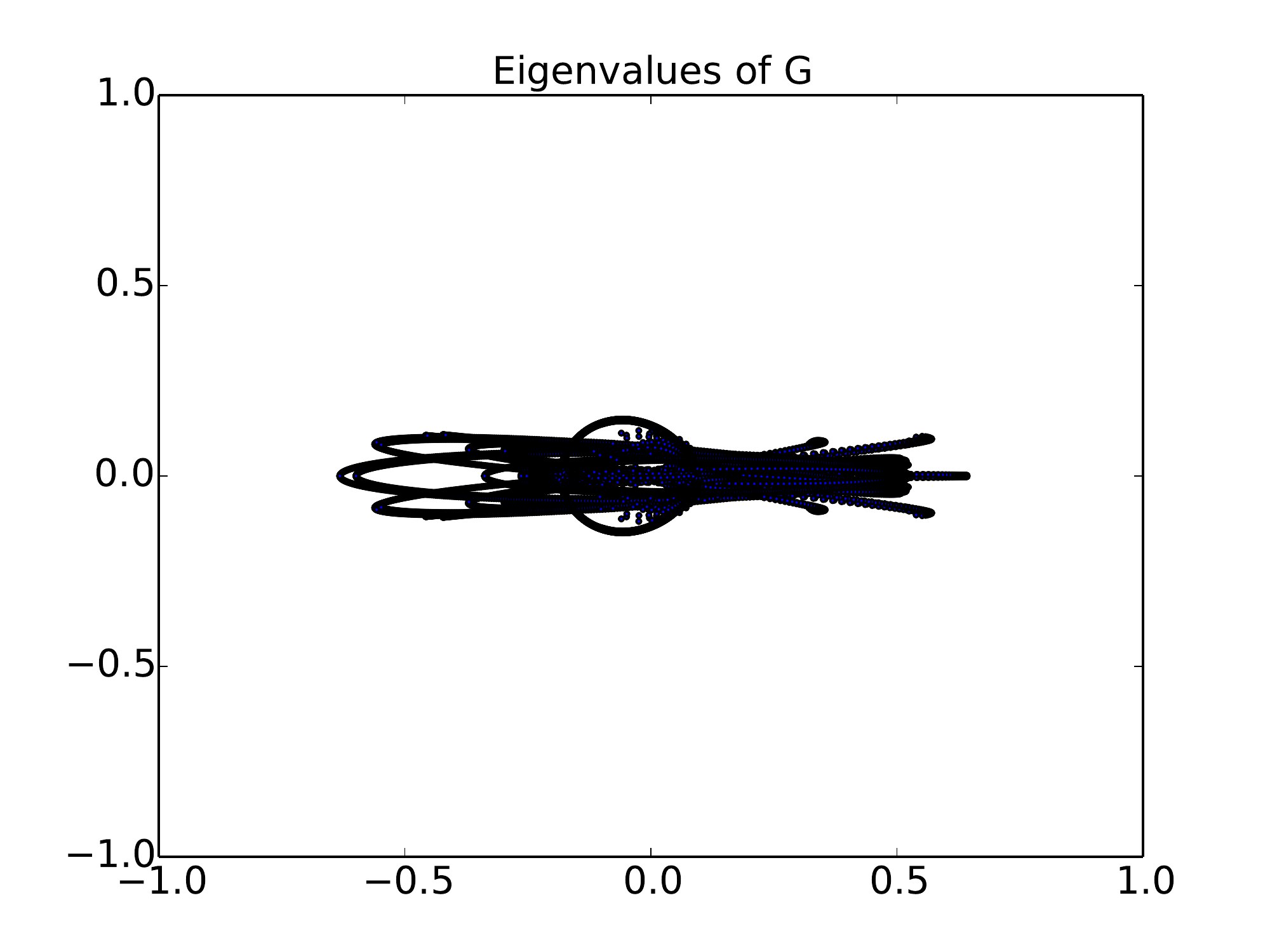}
\includegraphics[width=0.49\textwidth]{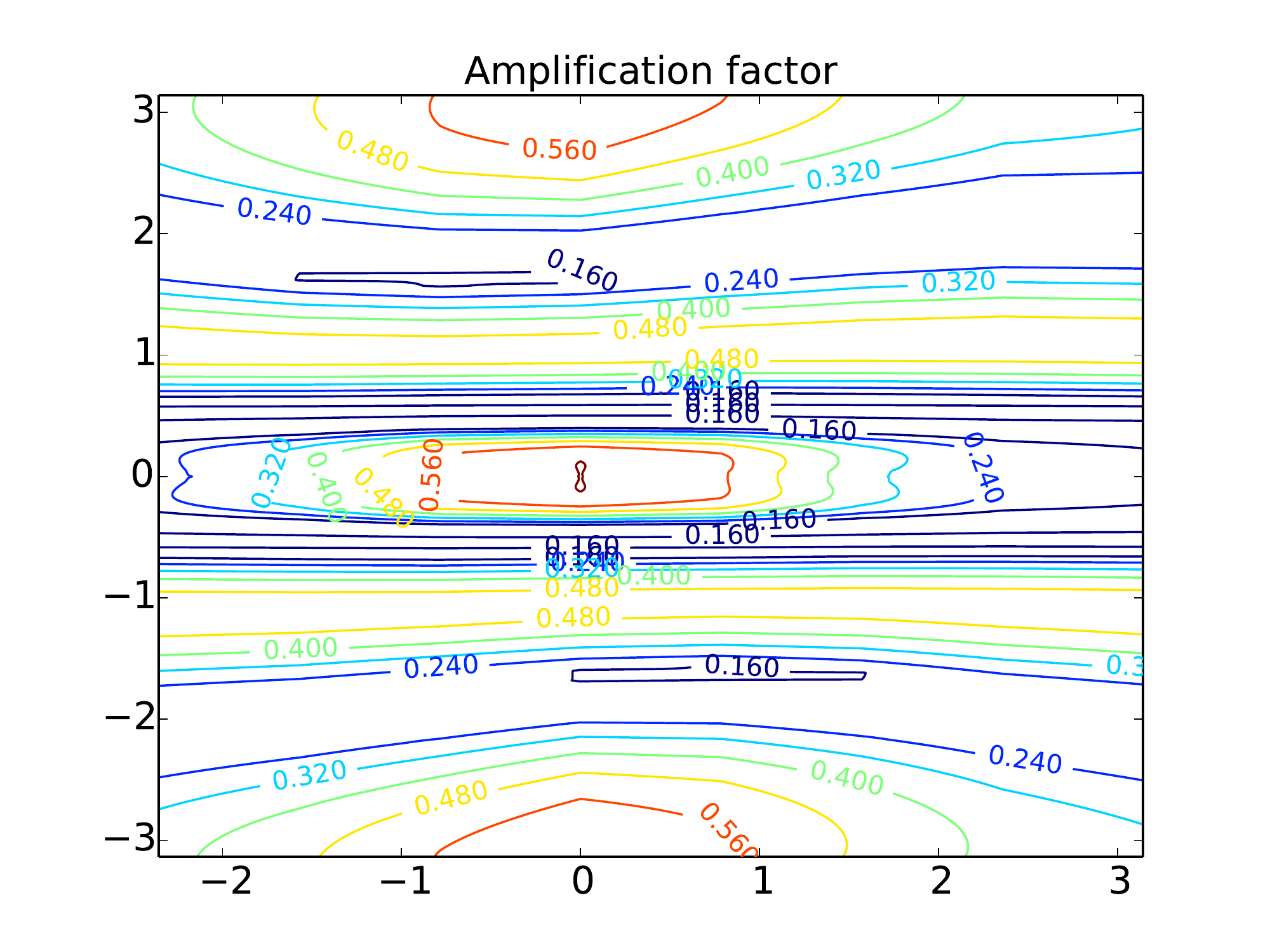}
\caption{\label{ampfactor-precARK}Spectrum and amplification factors for different wavenumbers, Mach 0.8, $\alpha = 0^{\circ}$, ARK3J with SGS preconditioner, $c$=200, $d$=0.5. Top: AR=1, $c^*$=14; Bottom: AR=100, $c^*$=4300.}
\end{figure}

\begin{table}
\centering
\begin{tabular}{l|ccc|ccc}
& \multicolumn{3}{c|}{SGS precond.} & \multicolumn{3}{c}{exact precond.} \\ \hline
$d$\qquad AR & 1 & 100 & 10000 & 1 & 100 & 10000\\ \hline
0.0 & 1 & 80 & 10500 & 1 & 80 & 8500 \\ \hline
0.1 & 6 & 900 & 95000 & 8 & 850 & 80000 \\ \hline
0.25 & 15 & 2200 & 220000 & 16 & 2050 & 200000 \\ \hline
0.5 & 30 & 4400 & 440000 & 27 & 4000 & 400000 \\ \hline
1.0 & 55 & 6500 & 800000 & 59 & 6800 & 790000 \\ \hline
\end{tabular}
\caption{\label{Prec-ARK-stability}Maximal $c^*$ for ARK3J, $c$=200, $M=0.8$, $\alpha = 0^{\circ}$.}
\end{table}

\begin{table}
\centering
\begin{tabular}{c|c|ccc|ccc}
& & \multicolumn{3}{c|}{SGS precond.} & \multicolumn{3}{c}{exact precond.} \\ \hline
$d$ & AR                  & 1 & 100 & 10000 & 1 & 100 & 10000\\ \hline
0.0 & $c^*$ opt      & 1 & 80 & 10500 & 1 & 80 & 8500 \\
    & $\rho({\bf M})$ opt & 0.9974 & 0.9708 & 0.9455 & 0.9480 & 0.9587 & 0.9557 \\
    & Sm. fct. opt        & 0.9439 & 0.9708 & 0.9455 & 0.9446 & 0.9587 & 0.9557 \\ \hline
0.1 & $c^*$ opt      & 5 & 900 & 95000 & 8 & 800 & 80000 \\
    & $\rho({\bf M})$ opt & 0.9878 & 0.7730 & 0.7739 & 0.6395 & 0.6434 & 0.6424 \\
    & Sm. fct. opt        & 0.7370 & 0.7730 & 0.7739 & 0.6533 & 0.6434 & 0.6424 \\ \hline
0.25 & $c^*$ opt      & 10 & 2100 & 220000 & 9 & 1900 & 190000 \\
    & $\rho({\bf M})$ opt & 0.9774 & 0.6932 & 0.6649 & 0.6028 & 0.5273 & 0.4282 \\
    & Sm. fct. opt        & 0.5480 & 0.6932 & 0.6649 & 0.5738 & 0.5273 & 0.4282 \\ \hline
0.5 & $c^*$ opt      & 14 & 4300 & 440000 & 14 & 2200 & 370000 \\
    & $\rho({\bf M})$ opt & 0.9722 & 0.6422 & 0.7411 & 0.4504 & 0.3158 & 0.3653 \\
    & Sm. fct. opt        & 0.4240 & 0.6422 & 0.7411 & 0.4504 & 0.3158 & 0.3653 \\ \hline
1.0 & $c^*$ opt      & 20 & 6100 & 700000 & 30 & 6300 & 650000 \\
    & $\rho({\bf M})$ opt & 0.9661 & 0.7122 & 0.6707 & 0.2653 & 0.5088 & 0.3984 \\
    & Sm. fct. opt        & 0.2908 & 0.7122 & 0.6707 & 0.2653 & 0.5088 & 0.3984 \\ \hline
\end{tabular}
\caption{\label{Prec-ARK-rates}Amplification and smoothing factors of ARK3J, $c$=200, $M=0.8$, $\alpha = 0^{\circ}$.}
\end{table}

We now consider preconditioned ARK3J with SGS and exact preconditioning. The Mach number is set to 0.8 and the angle of attack to zero degrees, which is the most difficult test case of the ones considered. Even so, it is possible to achieve convergence at all aspect ratios with a large physical CFL $c$=200. With regards to stability, we show the maximal possible $c^*$ in table~\ref{Prec-ARK-stability}. We can see that this is dramatically improved compared to the unpreconditioned method, but it remains finite, as predicted by the theory. We furthermore notice that the choice of $d$ in the cutoff function \eqref{cutoff} is important. In particular, the smaller we choose $d$, meaning the smaller we allow eigenvalues to be, the less stable the method will be. Maximal $c^*$ is approximately proportional to the aspect ratio and to $d$. The eigenvalues and contours of smoothing factor for $d$=0.5 are also illustrated in Figure~\ref{ampfactor-precARK} for aspect ratios 1 and 100, respectively. Clustering of the eigenvalues along the real axis is observed indicating good convergence.

For each value of $d$ considered, $c^*$ was optimised ($c^*$ opt) to minimise the smoothing factor (SM fct. opt). The results are shown in table~\ref{Prec-ARK-rates}. Optimal smoothing factors improve as $d$ is increased. Preconditioning with the exact inverse affords better smoothing factors than SGS preconditioning. With SGS preconditioning optimal smoothing factors at AR=1 are on the whole lower at than those at large AR. Conversely smoothing factors at AR=1 are equal to or higher than those at large AR when exact preconditioning is used. In general, optimal $c^*$ is close to maximal $c^*$. 

\subsubsection{Additive W methods}

 \begin{figure}[h]
\centering
\includegraphics[width=0.49\textwidth]{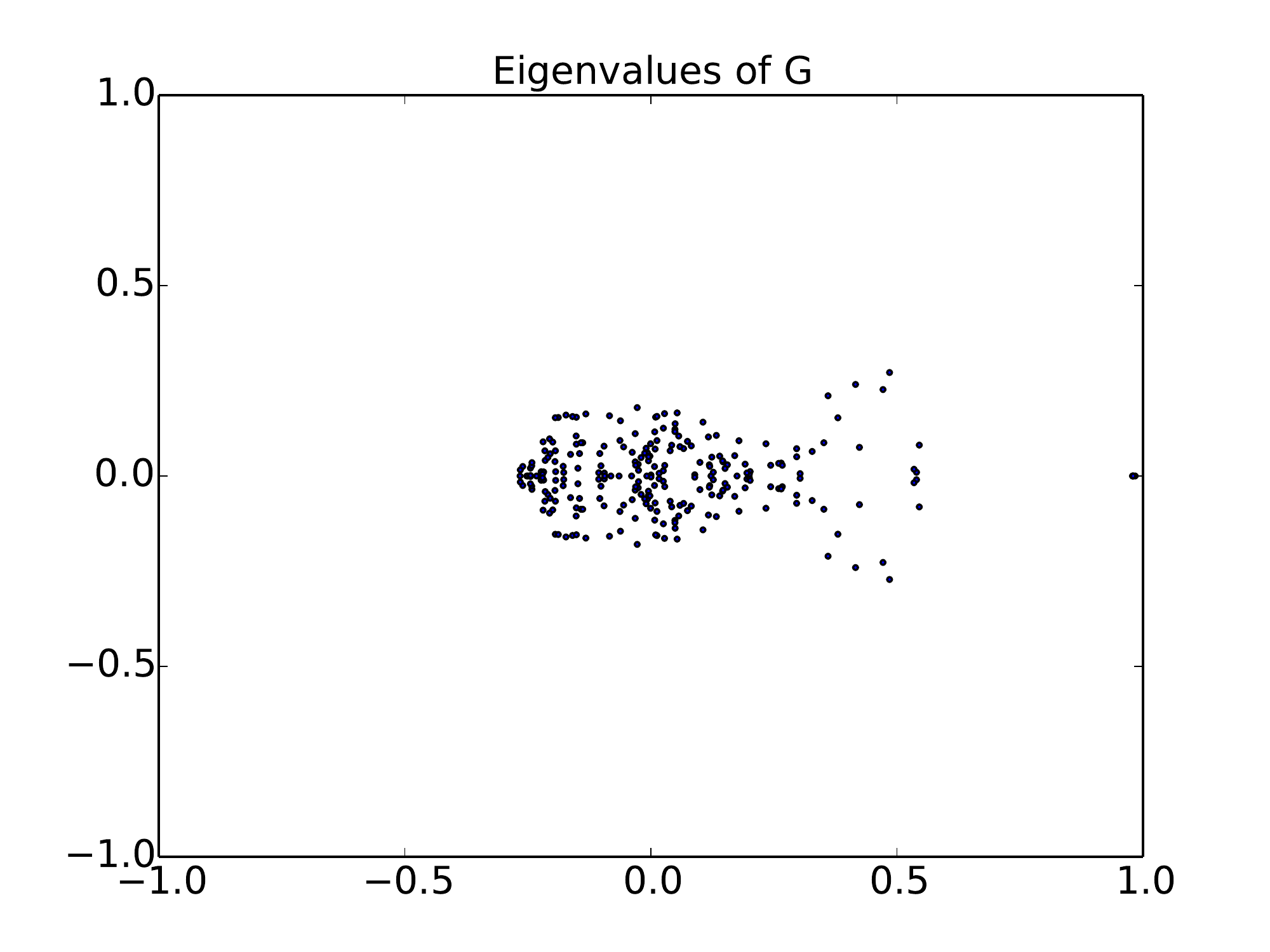}
\includegraphics[width=0.49\textwidth]{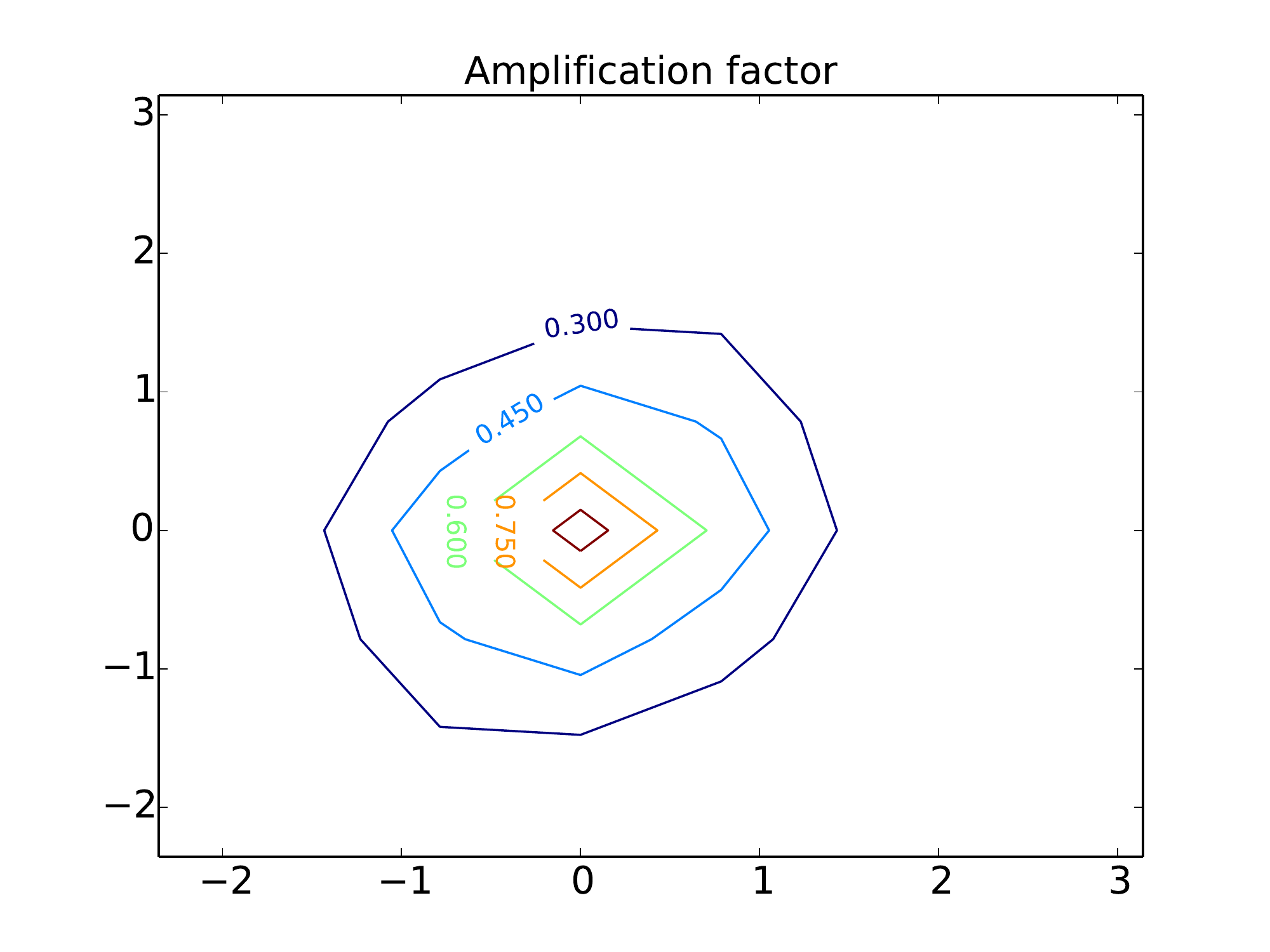}
\includegraphics[width=0.49\textwidth]{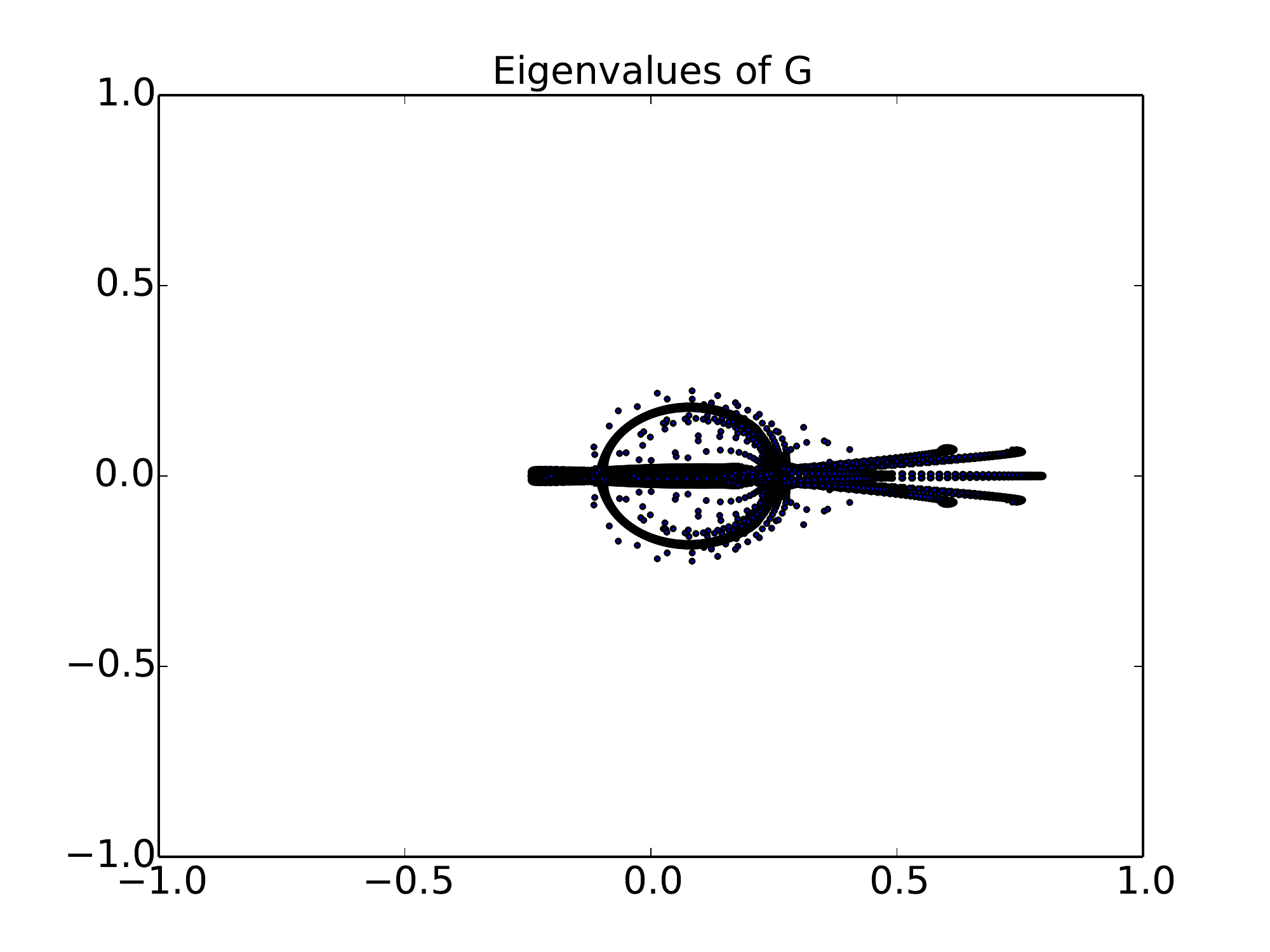}
\includegraphics[width=0.49\textwidth]{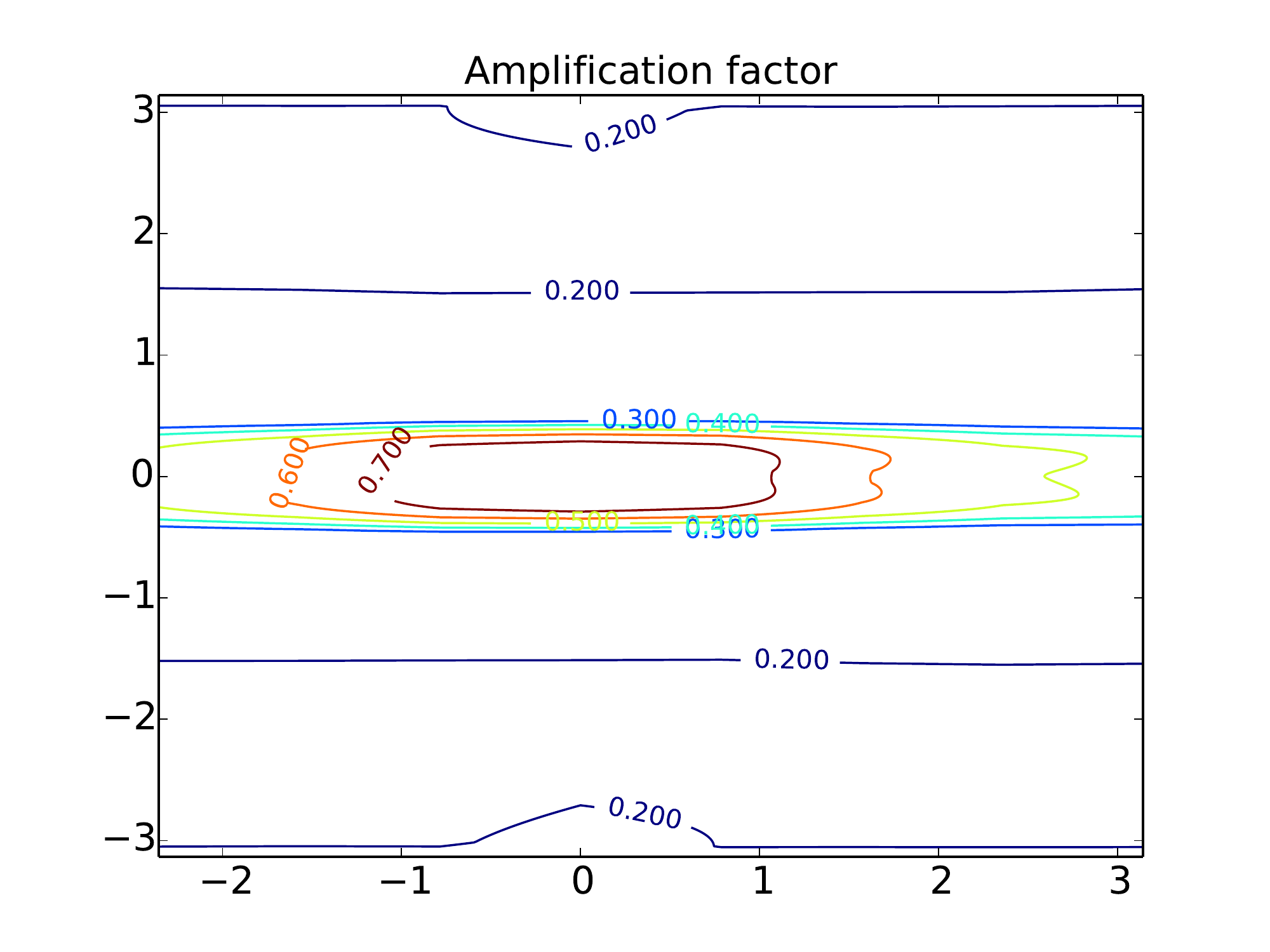}
\caption{\label{ampfactor-AW3}Spectrum and amplification factors for different wavenumbers, Mach 0.8, $\alpha=0^{\circ}$, AW3 with SGS preconditioner, $c$=200, $d$=0.5, $\eta$=0.8. Top: AR=1, $c^*$=3; Bottom: AR=100, $c^*=1e6$.}
\end{figure}

\begin{table}
\centering
\begin{tabular}{c|ccc|ccc}
& \multicolumn{3}{c|}{SGS precond.} & \multicolumn{3}{c}{exact precond.} \\ \hline
$d$\ AR & 1 & 100 & 10000 & 1 & 100 & 10000\\ \hline
0.0 & 8 & 8 & 8 & 8 & 8 & 8 \\ \hline
0.1 & 13 & 13 & 13 & 11 & 13 & 13 \\ \hline
0.25 & 30 & 2100 & A & 47 & 98 & 98 \\ \hline
0.5 & A & A & A & A & A & A \\ \hline
1.0 & A & A & A & A & A & A \\ \hline
\end{tabular}
\caption{\label{AW-stability}Maximal $c^*$ for AW3, $\eta$=0.8, $c$=200, $M=0.8$, $\alpha = 0^{\circ}$. A implies that no bound was observed.}
\end{table}

\begin{table}
\centering
\begin{tabular}{c|c|ccc|ccc}
& & \multicolumn{3}{c|}{SGS precond.} & \multicolumn{3}{c}{exact precond.} \\ \hline
$d$ & AR                  & 1 & 100 & 10000 & 1 & 100 & 10000\\ \hline
0.0 & $c^*$ opt      & 3 & 8 & 8 & 3 & 8 & 8 \\
    & $\rho({\bf M})$ opt & 0.9836 & 0.9452 & 0.9441 & 0.9781 & 0.9440 & 0.9440 \\
    & Sm. fct. opt        & 0.3046 & 0.9452 & 0.9441 & 0.3046 & 0.9440 & 0.9440 \\ \hline
0.1 & $c^*$ opt      & 3 & 12 & 13 & 3 & 12 & 12 \\
    & $\rho({\bf M})$ opt & 0.9837 & 0.9217 & 0.9140 & 0.9781 & 0.9188 & 0.9188 \\
    & Sm. fct. opt        & 0.2969 & 0.9217 & 0.9140 & 0.2933 & 0.9188 & 0.9188 \\ \hline
0.25 & $c^*$ opt      & 3 & 240 & 500 & 3 & 70 & 70 \\
    & $\rho({\bf M})$ opt & 0.9838 & 0.7157 & 0.6813 & 0.9781 & 0.6831 & 0.6831 \\
    & Sm. fct. opt        & 0.2818 & 0.7157 & 0.6813 & 0.2787 & 0.6831 & 0.6831 \\ \hline
0.5 & $c^*$ opt      & 3 & $>1e6$ & $>1e6$ & 4 & $>1e6$ & $>1e6$ \\
    & $\rho({\bf M})$ opt & 0.9841 & 0.7965 & 0.7903 & 0.9740 & 0.3473 & 0.3473 \\
    & Sm. fct. opt        & 0.2686 & 0.7965 & 0.7903 & 0.2669 & 0.3473 & 0.3473 \\ \hline
1.0 & $c^*$ opt      & 7 & $>1e6$ & $>1e6$ & 7 & $>1e6$ & $>1e6$ \\
    & $\rho({\bf M})$ opt & 0.9765 & 0.8845 & 0.8822 & 0.9506 & 0.4045 & 0.4045 \\
    & Sm. fct. opt        & 0.2564 & 0.8845 & 0.8822 & 0.2627 & 0.4045 & 0.4045 \\ \hline
\end{tabular}
\caption{\label{AW}Amplification and smoothing factors of AW3, $\eta$=0.8, $c$=200, $M=0.8$, $\alpha = 0^{\circ}$.}
\end{table}

Results for AW3 with SGS and exact preconditioning are shown in Tables~\ref{AW-stability} and \ref{AW}. Again, the Mach number is set to 0.8, the physical CFL $c$=200 and the angle of attack to zero degrees. We set $\eta=0.8$. An $A$ means that no bound on $c^*$ was observed. As we can see, as long as $d$ is chosen sufficiently large, the methods are practically A-stable, as suggested by the theory. Surprisingly, for $d$ small, stability is worse than for the preconditioned ARK methods. 
This is also illustrated in Figure~\ref{ampfactor-AW3} for for Mach 0.5 and aspect ratios 1 and 100, respectively. As with ARK3J, the eigenvalues are clustered along the real axis.

A slightly more complex picture emerges when the optimal smoothing factor is considered.
At AR=1, the AW3 scheme attains very low optimal smoothing factors of around 0.3 at all values of $d$ while the ARK3J scheme smoothing factors improved with increasing $d$.
Comparing SGS preconditioning in both schemes, the optimal smoothing factors obtained by AW3 are slightly higher than ARK3J.
Using exact preconditioning in both schemes at AR=100 and 10000, AW3 and ARK3J obtain comparable smoothing factors.
Regarding the optimal $c^*$, it is generally lower than with ARK3J except for $d \geq 0.5$ and AR$>1$.

\subsubsection{Comparison of AW schemes and choice of $\eta$}

One important question is the optimal choice of the additional parameter $\eta$ in the W methods. 
Based on the AW3 results in Table \ref{AW-stability} it was decided to focus on two values of $d$: $d$=0.1 where limited stability was observed, and $d$=0.5 where $A$-stability was observed.
Only SGS preconditioning was used.
For each W scheme and value of $d$, optimal values of $c^*$, $\eta$ and amplification and smoothing factors were determined.
These are presented in Table \ref{AWa0} for initial conditions Mach=0.8, $\alpha=0^{\circ}$ and in Table \ref{AWa45} for initial conditions Mach=0.8, $\alpha=45^{\circ}$.

Looking just at Table \ref{AWa0}, the optimal value of $\eta$ is low, either 0.4 or 0.5 (with one case of 0.7), when $d$=0.5.
When $d$=0.1, the optimal $\eta$ depends on AR: for AR=1, optimal values of $\eta$ are 0.5 or 0.6 and for AR=100 and 10000 the values are higher, mostly 0.8.
Looking at Table \ref{AWa45}, the optimal value of $\eta$ is independent of $d$ and the choice of scheme but not of AR.
The optimal value of $\eta$ appears to be somewhat dependent on the initial conditions and other free parameters but independent of the specific W scheme.
Furthermore, the optimisation process demonstrated (not all results are shown for brevity) that the W schemes are all stable within a range: $0.5 \lessapprox \eta \lessapprox 0.9$ but the maximal $c^*$ varies with $\eta$ within the range.
As shown in Table \ref{AW}, fixing $\eta=0.8$ across all tests results in a stable but sub-optimal scheme.
Looking at the relative performance of different W schemes in Tables \ref{AWa0} and \ref{AWa45}, it is apparent that they all obtain similar optimal smoothing factors at similar $c^*$ values.
Therefore, AW3 is the best scheme as it uses only three stages.

The discrete Fourier analysis suggests that the preconditioned ARK3J and additive W schemes should theoretically achieve very good smoothing factors under challenging flow conditions and on high aspect ratio grids.
Moreover, in the W schemes the eigenvalue limiting parameter $d$ plays an important role: for $d \geq 0.5$ and AR$>1$ the allowable $c^*$ is unlimited, while for smaller $d$ or AR=1 the optimal $c^*$ is finite and smaller than that found for preconditioned ARK3J.

\begin{table}
\centering
\begin{tabular}{c|c|ccc|ccc}
& scheme & \multicolumn{3}{c|}{AW3} & \multicolumn{3}{c}{AW51} \\ \hline
$d$ & AR                  & 1 & 100 & 10000 & 1 & 100 & 10000 \\ \hline
0.1 & $c^*$ opt           & 3 & 12 & 13 & 3 & 8 & 11 \\
    & $\eta$ opt          & 0.6 & 0.8 & 0.8 & 0.6 & 0.8 & 0.9 \\
    & $\rho({\bf M})$ opt & 0.9823 & 0.9217 & 0.9140 & 0.9823 & 0.9454 & 0.9824 \\
    & Sm. fct. opt        & 0.2604 & 0.9217 & 0.9140 & 0.2624 & 0.9454 & 0.9824 \\ \hline
0.5 & $c^*$ opt           & 3 & $>1e6$ & $>1e6$ & 4 & 30 & $>1e6$ \\
    & $\eta$ opt          & 0.4 & 0.5 & 0.5 & 0.5 & 0.5 & 0.7 \\
    & $\rho({\bf M})$ opt & 0.9809 & 0.6918 & 0.6830 & 0.9776 & 0.8598 & 0.7656 \\
    & Sm. fct. opt        & 0.2630 & 0.6918 & 0.6830 & 0.2611 & 0.8598 & 0.7656 \\ \hline
& scheme & \multicolumn{3}{c|}{AW52} & \multicolumn{3}{c}{AW5J} \\ \hline
$d$ & AR                  & 1 & 100 & 10000 & 1 & 100 & 10000 \\ \hline
0.1 & $c^*$ opt           & 3 & 8 & 10 & 3 & 9 & 10 \\
    & $\eta$ opt          & 0.6 & 0.8 & 0.8 & 0.5 & 0.7 & 0.8 \\
    & $\rho({\bf M})$ opt & 0.9823 & 0.9456 & 0.9323 & 0.9815 & 0.9389 & 0.9323 \\
    & Sm. fct. opt        & 0.2256 & 0.9456 & 0.9323 & 0.2064 & 0.9389 & 0.9323 \\ \hline
0.5 & $c^*$ opt           & 3 & $>1e6$ & $>1e6$ & 3 & $>1e6$ & $>1e6$ \\
    & $\eta$ opt          & 0.5 & 0.5 & 0.5 & 0.4 & 0.5 & 0.5 \\
    & $\rho({\bf M})$ opt & 0.9818 & 0.7016 & 0.6934 & 0.9809 & 0.6991 & 0.6907 \\
    & Sm. fct. opt        & 0.1762 & 0.7016 & 0.6934 & 0.1721 & 0.6991 & 0.6907 \\ \hline
\end{tabular}
\caption{\label{AWa0}Optimal $\eta$, $c^*$, amplification and smoothing factors of all AW schemes, $c$=200, $M=0.8$, $\alpha = 0^{\circ}$, SGS preconditioning.}
\end{table}

\begin{table}
\centering
\begin{tabular}{c|c|ccc|ccc}
& scheme & \multicolumn{3}{c|}{AW3} & \multicolumn{3}{c}{AW51} \\ \hline
$d$ & AR                  & 1 & 100 & 10000 & 1 & 100 & 10000 \\ \hline
0.1 & $c^*$ opt           & 3 & 900 & 400 & 3 & 1200 & 300 \\
    & $\eta$ opt          & 0.6 & 0.8 & 0.9 & 0.6 & 0.8 & 0.8 \\
    & $\rho({\bf M})$ opt & 0.9804 & 0.4481 & 0.4367 & 0.9804 & 0.4545 & 0.4390 \\
    & Sm. fct. opt        & 0.2642 & 0.4481 & 0.4367 & 0.2654 & 0.4545 & 0.4390 \\ \hline
0.5 & $c^*$ opt           & 3 & $>1e6$ & $>1e6$ & 3 & $>1e6$ & $>1e6$ \\
    & $\eta$ opt          & 0.6 & 0.8 & 0.9 & 0.6 & 0.8 & 0.9 \\
    & $\rho({\bf M})$ opt & 0.9809 & 0.4441 & 0.4363 & 0.9804 & 0.4484 & 0.4351 \\
    & Sm. fct. opt        & 0.2642 & 0.4441 & 0.4363 & 0.2654 & 0.4484 & 0.4351 \\ \hline
& scheme & \multicolumn{3}{c|}{AW52} & \multicolumn{3}{c}{AW5J} \\ \hline
$d$ & AR                  & 1 & 100 & 10000 & 1 & 100 & 10000 \\ \hline
0.1 & $c^*$ opt           & 3 & $>1e6$ & 500 & 3 & 300 & 200 \\
    & $\eta$ opt          & 0.5 & 0.8 & 0.8 & 0.5 & 0.8 & 0.9 \\
    & $\rho({\bf M})$ opt & 0.9798 & 0.4378 & 0.4107 & 0.9798 & 0.5460 & 0.5220 \\
    & Sm. fct. opt        & 0.1750 & 0.4378 & 0.4107 & 0.1526 & 0.5460 & 0.5220 \\ \hline
0.5 & $c^*$ opt           & 3 & $>1e6$ & $>1e6$ & 3 & 1100 & 800 \\
    & $\eta$ opt          & 0.5 & 0.8 & 0.9 & 0.5 & 0.8 & 0.9 \\
    & $\rho({\bf M})$ opt & 0.9799 & 0.4710 & 0.3957 & 0.9799 & 0.5427 & 0.4205 \\
    & Sm. fct. opt        & 0.1750 & 0.4710 & 0.3957 & 0.1527 & 0.5427 & 0.4205 \\ \hline
\end{tabular}
\caption{\label{AWa45}Optimal $\eta$, $c^*$, amplification and smoothing factors of all AW schemes, $c$=200, $M=0.8$, $\alpha = 45^{\circ}$, SGS preconditioning.}
\end{table}

\section{Numerical results}

We now proceed to tests on the RANS equations and use a FAS scheme as the iterative solver. We employ the Fortran code uflo103 to compute flows around pitching airfoils. All computations are run on Ubuntu 16.04 on a single core of an 8-core Intel i7-3770 CPU at 3.40GHz with 8 GB of memory.

C-type grids are employed, where the half of the cells that are closer to the boundary in $y$-direction get a special boundary layer scaling. To obtain initial conditions for the unsteady simulation, far field values are used from which a steady state is computed. The first unsteady time step does not use BDF-2, but implicit Euler as a startup for the multistep method. From then on, BDF-2 is employed. We look at the startup phase to evaluate the performance of steady state computations and at the second overall timestep, meaning the first BDF-2 step, to evaluate performance for the unsteady case. 

\begin{figure}[h]
\centering
\includegraphics[height=5.0cm]{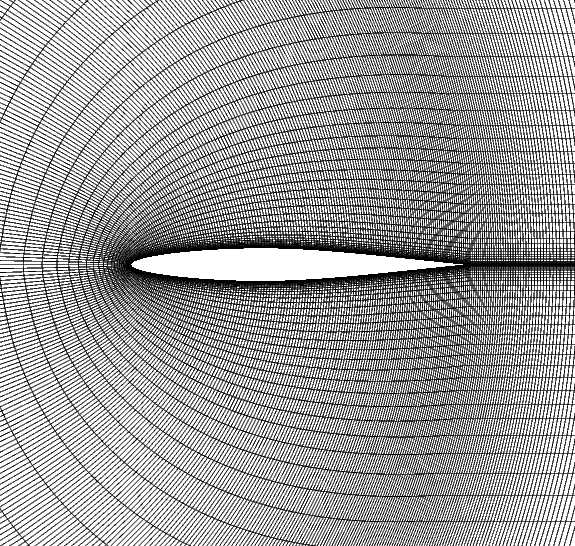}
\includegraphics[height=5.0cm]{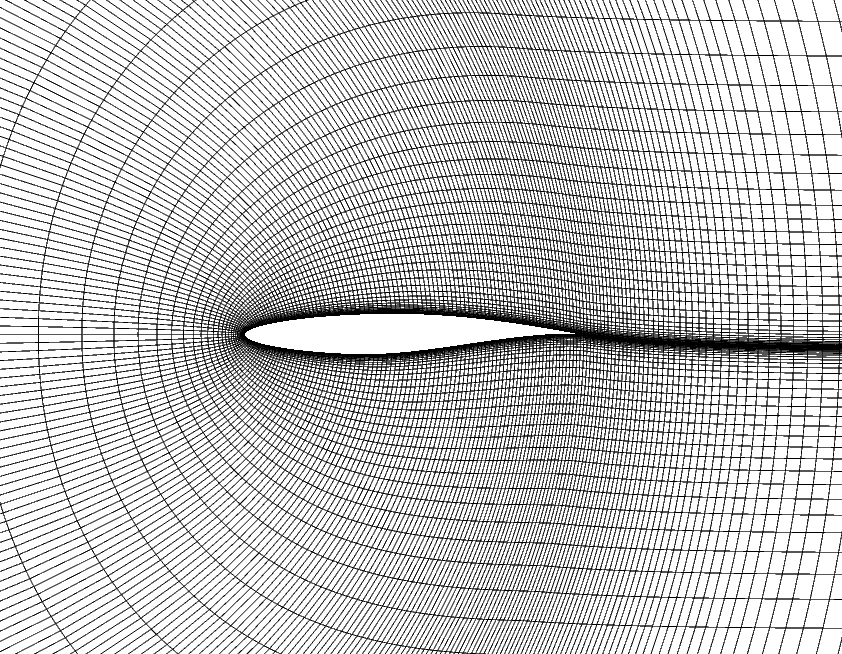}
\caption{\label{grids}Zoom of grids around NACA 64A010 and RAE 2822 airfoils.}
\end{figure}

As a first test case, we consider the flow around the NACA 64A010 pitching and plunging airfoil at a Mach number 0.796. The grid is illustrated in Figure~\ref{grids}. For the pitching, we use a frequency of 0.202 and an amplitude of $1.01^{\circ}$. 36 timesteps per cycle (pstep) are chosen. The Reynolds number is $10^6$ and the Prandtl number is 0.75. The grid is a C-mesh with $512\times 64$ cells and maximum aspect ratio of $6.31e6$. As a second test case, we look at the pitching RAE 2822 airfoil at a Mach number of 0.75. The grid is illustrated in Figure~\ref{grids}. For the pitching, we use a frequency of 0.202 and an amplitude of $1.01^{\circ}$ and pstep=36. The grid has $320\times 64$ cells and maximum aspect ratio of $8.22e6$.

The results of the Fourier analysis suggest that the most interesting schemes are SGS preconditioned ARK3J and the various AW schemes. A first thing to note is that due to nonlinear effects, the schemes need to be tweaked from the linear to the nonlinear case. In particular, it is necessary to start with a reduced pseudo CFL number $c^*$. We restrict it to 20 for the first two iterations.

\subsection{Choice of parameters}

\begin{table}
\centering
\begin{tabular}{c|cc|cc}
& \multicolumn{2}{c|}{Steady} & \multicolumn{2}{c}{Unsteady} \\ \hline
$\eta$ & $c^*$ & Conv. rate & $c^*$ & Conv. rate \\ \hline
0.4 & 10 & 0.8774 & 10 & 0.8770 \\
0.5 & 10000 & 0.8616 & 10000 & 0.8071 \\
0.6 & 10000 & 0.8629 & 10000 & 0.8183 \\
1.0 & 10000 & 0.8759 & 10000 & 0.8546
\end{tabular}
\caption{\label{uflo-rates-naca-d05}Maximal $c^*$ and convergence rates of UFLO103 for NACA 64A010 test case, $d=0.5$}
\end{table}

\begin{table}
\centering
\begin{tabular}{c|cc|cc}
& \multicolumn{2}{c|}{Steady} & \multicolumn{2}{c}{Unsteady} \\ \hline
$\eta$ & $c^*$ & Conv. rate & $c^*$ & Conv. rate \\ \hline
0.5 & 10000 & 0.8561 & 100 & 0.7517 (90) \\
0.6 & 10000 & 0.8583 & 100 & 0.7732
\end{tabular}
\caption{\label{uflo-rates-naca-d01}Maximal $c^*$ and convergence rates of UFLO103 for NACA 64A010 test case, $d=0.1$}
\end{table}

\begin{table}
\centering
\begin{tabular}{c|cc|cc}
& \multicolumn{2}{c|}{Steady} & \multicolumn{2}{c}{Unsteady} \\ \hline
$\eta$ & $c^*$ & Conv. rate & $c^*$ & Conv. rate \\ \hline
0.4 & 10 & 0.8787 & 10 & 0.8828 \\
0.5 & 10000 & 0.8506 & 100 & 0.7503 (90) \\
0.6 & 10000 & 0.8533 & 100 & 0.7721
\end{tabular}
\caption{\label{uflo-rates-naca-d005}Maximal $c^*$ and convergence rates of UFLO103 for NACA 64A010 test case, $d=0.05$}
\end{table}

\begin{table}
\centering
\begin{tabular}{c|cc|cc}
& \multicolumn{2}{c|}{Steady} & \multicolumn{2}{c}{Unsteady} \\ \hline
$\eta$ & $c^*$ & Conv. rate & $c^*$ & Conv. rate \\ \hline
0.4 & 10000 & 0.8567 & 10000 & 0.8228 \\
0.5 & 10000 & 0.8581 & 10000 & 0.8424 \\
0.6 & 10000 & 0.8631 & 10000 & 0.8574 \\
0.7 & 10000 & 0.8691 & 10000 & 0.8688 \\
0.8 & 10000 & 0.8740 & 10000 & 0.8766 \\
0.9 & 10000 & 0.8744 & 10000 & 0.8817 \\
1.0 & 10000 & 0.8804 & 10000 & 0.8866
\end{tabular}
\caption{\label{uflo-rates-rae-d05}Maximal $c^*$ and convergence rates of UFLO103 for RAE 2822 test case, $d=0.5$}
\end{table}

\begin{table}
\centering
\begin{tabular}{c|cc|cc}
& \multicolumn{2}{c|}{Steady} & \multicolumn{2}{c}{Unsteady} \\ \hline
$\eta$ & $c^*$ & Conv. rate & $c^*$ & Conv. rate \\ \hline
0.4 & 400 & 0.8429 & 60 & 0.7996 \\
0.5 & 10000 & 0.8359 & 70 & 0.8000 \\
0.8 & 10000 & 0.8471 & 60 & 0.8304
\end{tabular}
\caption{\label{uflo-rates-rae-d01}Maximal $c^*$ and convergence rates of UFLO103 for RAE 2822 test case, $d=0.1$}
\end{table}

\begin{table}
\centering
\begin{tabular}{c|cc|cc}
& \multicolumn{2}{c|}{Steady} & \multicolumn{2}{c}{Unsteady} \\ \hline
$\eta$ & $c^*$ & Conv. rate & $c^*$ & Conv. rate \\ \hline
0.4 & 100 & 0.8470 & 60 & 0.7895 \\
0.5 & 10000 & 0.8224 & 60 & 0.7972 \\
0.6 & 10000 & 0.8281 & 60 & 0.8045 \\
0.7 & 10000 & 0.8326 & 70 & 0.8029 \\
\end{tabular}
\caption{\label{uflo-rates-rae-d005}Maximal $c^*$ and convergence rates of UFLO103 for RAE 2822 test case, $d=0.05$}
\end{table}

In the AW methods, there are now three interdependent parameters to choose: $\eta$, $d$ and $c^*$, the CFL number in pseudo time. We start by fixing $d$. Choosing $d=0$ does not cause instability per se, but it leads to a stall in the iteration away from the solution. The convergence rates for $d=0.05$, $d=0.1$ and $d=0.5$ for the NACA and the RAE test case can be seen in tables~\ref{uflo-rates-naca-d05}-\ref{uflo-rates-rae-d005}. The largest $c^*$ tried is 10000 in all cases. If the number reported is smaller, it implies that it is the largest for which the methods are convergent. A number in parentheses e.g. (90) after the convergence rate means that the rate was calculated for the first 90 iterations, after which convergence stalled. Only stable values of $\eta$ are reported for brevity. The schemes are stable within a certain range, $0.5 \leq \eta \leq 0.9$, which tallies with the Fourier analysis results.

Qualitatively, we observe the following behavior:
\begin{itemize}
\item Increasing $d$ makes the schemes slower to converge and more stable
\item This effect is stronger for the unsteady system
\item If $\eta$ is too small, we get instability
\item Decreasing $\eta$ within the stable region will improve the convergence rate
\end{itemize}
 
We thus suggest two different modes of operation:
\begin{enumerate}
\item The robust mode: Choose $d=0.5$, $\eta=0.5$ and $c^*$ very large
\item The fast mode: Choose $d=0.05$, $\eta=0.5$ and $c^*$=100
\end{enumerate}
The robust mode trades some convergence rate for more robustness. 

The numerical experiments find somewhat different optimal values of $\eta$ to those found in the discrete Fourier analysis. Possible reasons for the discrepancies include the linearisations used in the discrete Fourier analysis and the non-cartesian meshes in the numerical experiments.

\subsection{Comparison of Schemes}

The linear analysis suggests that preconditioned ARK3J is competitive with the preconditioned W methods in terms of smoothing power. However, its application requires choosing $c^*$ within a stability limit whereas the W methods are $A$-stable for a certain range of $d$. To test the stability limit of the ARK schemes, we apply preconditioned ARK3J and ARK51 to the pitching NACA airfoil test case. The ARK3J method becomes unstable for $c^*>1$, whereas ARK51 can be run with $c^*$=3. However, both methods are completely uncompetitive with convergence rates of 0.999. Hereafter we compare only the AW schemes.

\begin{figure}[h]
\centering
\subfigure[Steady]{
\includegraphics[width=0.48\textwidth]{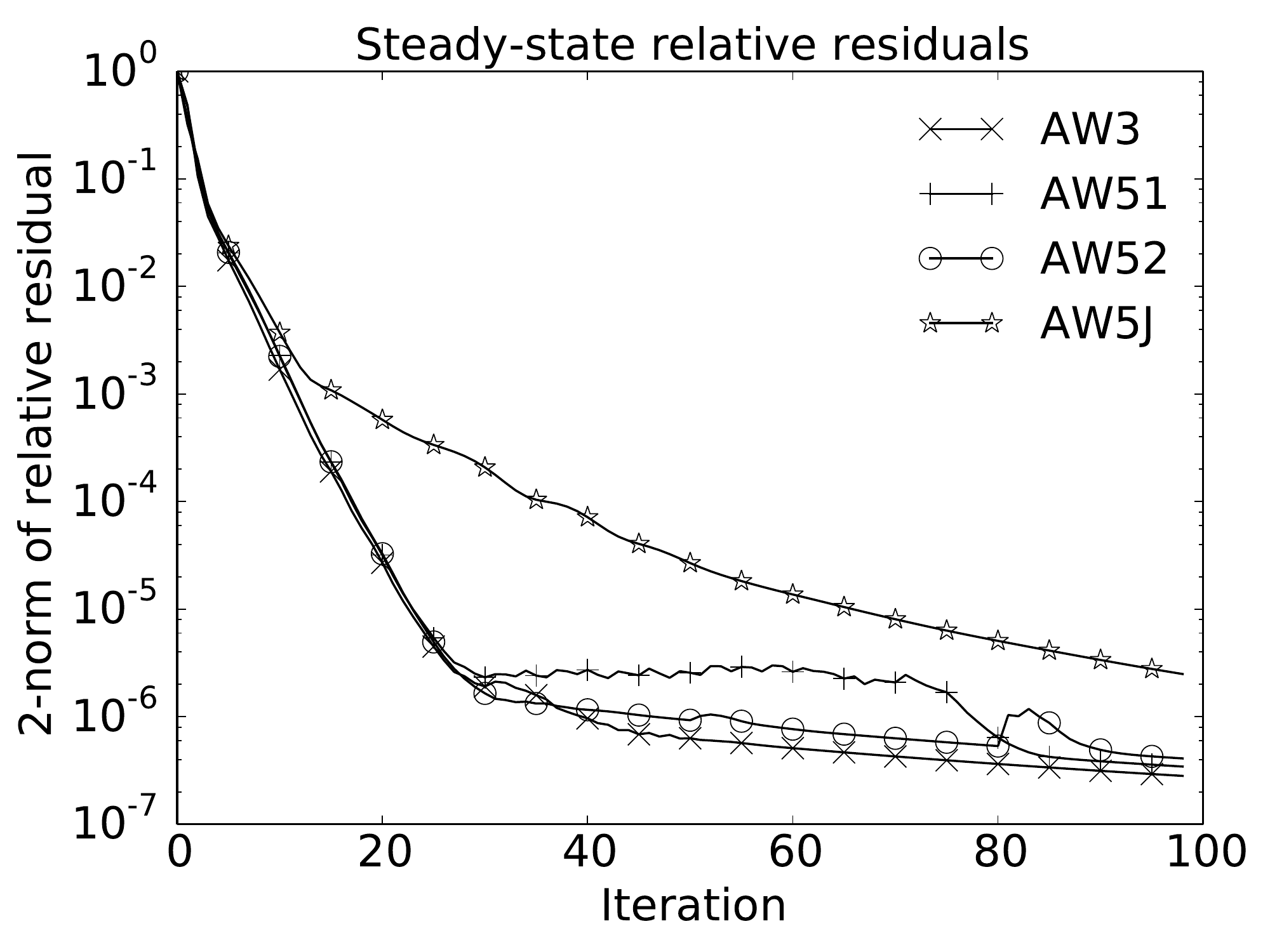}}
\subfigure[Unsteady]{
\includegraphics[width=0.48\textwidth]{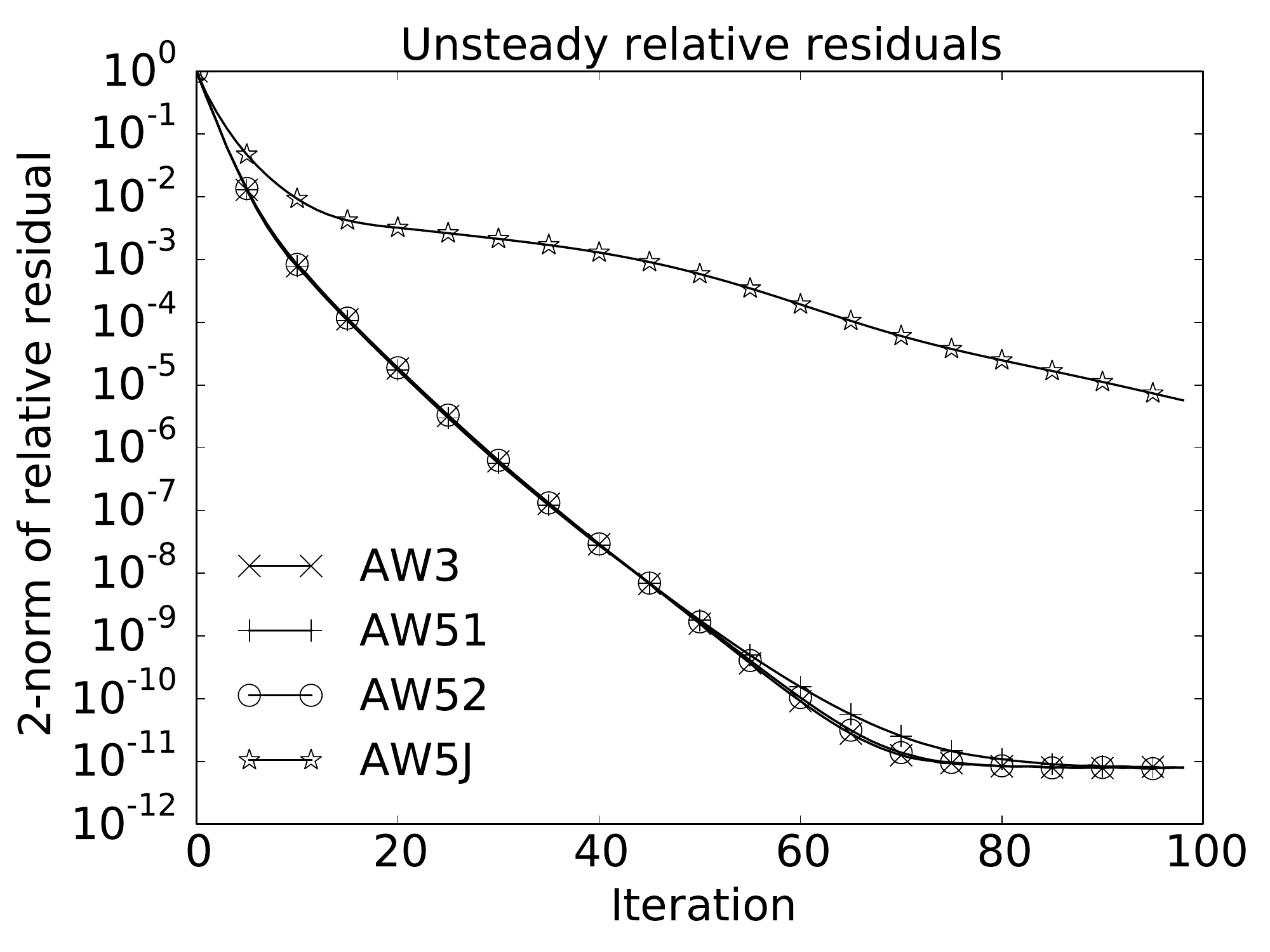}}
\caption{\label{conv_naca_all_d005}Convergence behavior for steady and unsteady flow around NACA64A010 airfoil for different SGS preconditioned W schemes, $d=0.05$.}
\end{figure}

\begin{figure}[h]
\centering
\subfigure[Steady]{
\includegraphics[width=0.48\textwidth]{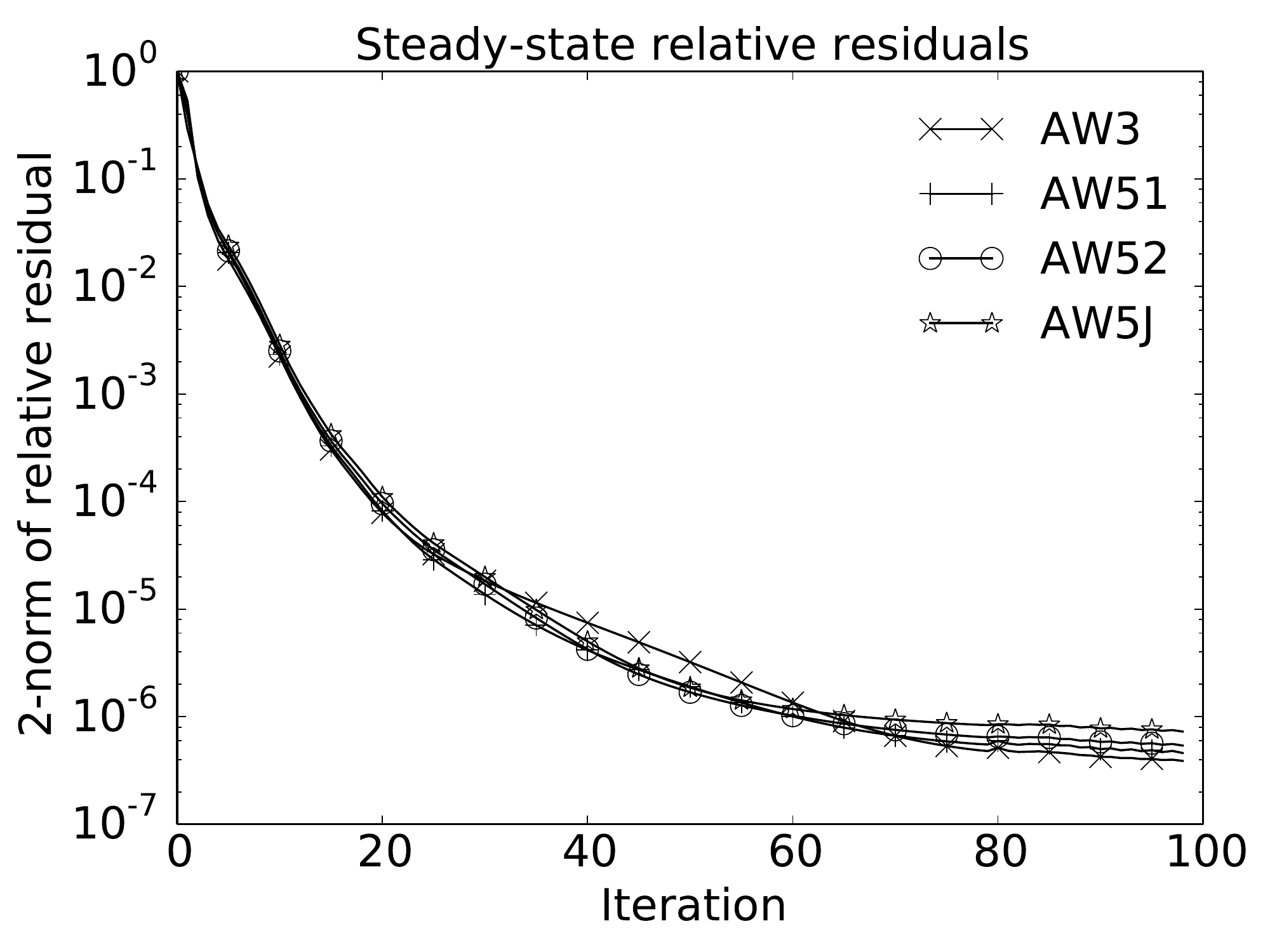}}
\subfigure[Unsteady]{
\includegraphics[width=0.48\textwidth]{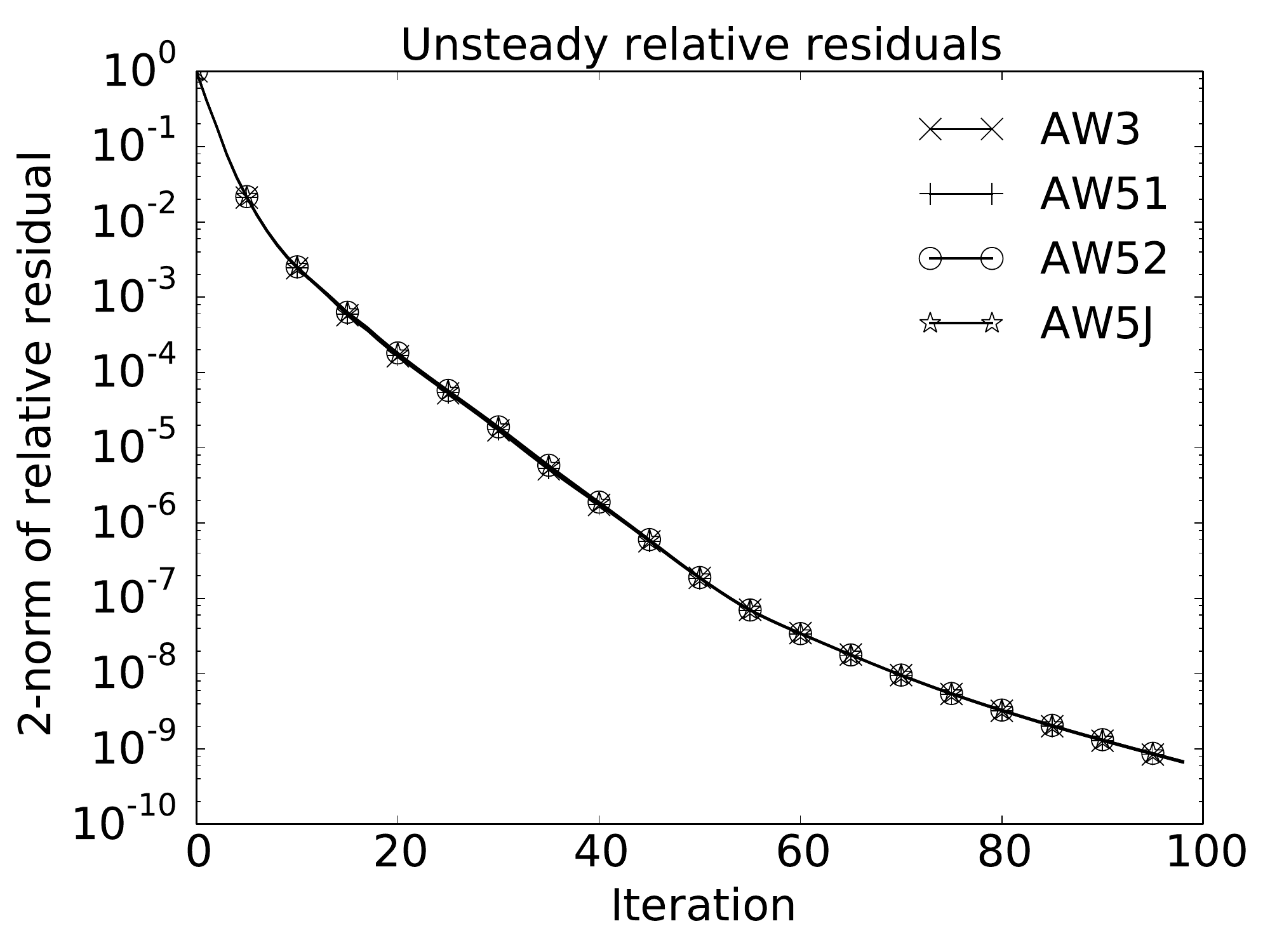}}
\caption{\label{conv_naca_all_d05}Convergence behavior for steady and unsteady flow around NACA64A010 airfoil for different SGS preconditioned W schemes, $d=0.5$.}
\end{figure}

\begin{table}
\centering
\begin{tabular}{cc|cc|cc}
 & & \multicolumn{2}{|c|}{$d=0.05$, $c^*$=100} & \multicolumn{2}{|c}{$d=0.5$, $c^*$=10000} \\ \hline
       &        & CPU[s] & av. conv. rate & CPU[s] & av. conv. rate \\ \hline
Steady & AW3   & 13.6 & 0.8586 & 14.5 & 0.8616 \\
       & AW5J  & 21.3 & 0.8775 & 21.4 & 0.8671 \\
       & AW51  & 21.0 & 0.8603 & 21.3 & 0.8632 \\ 
       & AW52  & 20.8 & 0.8618 & 21.5 & 0.8645 \\ \hline
Unsteady & AW3 & 19.0 & 0.7724 & 18.9 & 0.8071 \\
BDF-2 & AW5J   & 28.6 & 0.8844 & 29.0 & 0.8074 \\ 
      & AW51   & 28.7 & 0.7725 & 29.6 & 0.8071 \\ 
      & AW52   & 28.5 & 0.7724 & 29.2 & 0.8074 \\ \hline
\end{tabular}
\caption{\label{uflo-naca}Performance of AW3, AW51 and AW52 for the pitching NACA 64A010 airfoil, $0^{\circ}$ angle of attack, 100 steady/unsteady iterations.}
\end{table}

We first look at the NACA airfoil and compare AW3, AW51, AW52 and AW5J for the two modes of operation: $d$=0.05 and $c^*$=100 versus $d$=0.5 and $c^*$=10000. The relative residuals for the initial steady state computation and for the second unsteady time step are plotted in Figure~\ref{conv_naca_all_d005} for $d$=0.05 and $c^*$=100 and in Figure \ref{conv_naca_all_d05} for $d$=0.5 and $c^*$=10000. The convergence rates and CPU times are summarized in Table~\ref{uflo-naca}. Faster convergence is obtained with $d$=0.05 for all schemes except AW5J.

\begin{figure}[h]
\centering
\subfigure[steady]{
\includegraphics[width=0.48\textwidth]{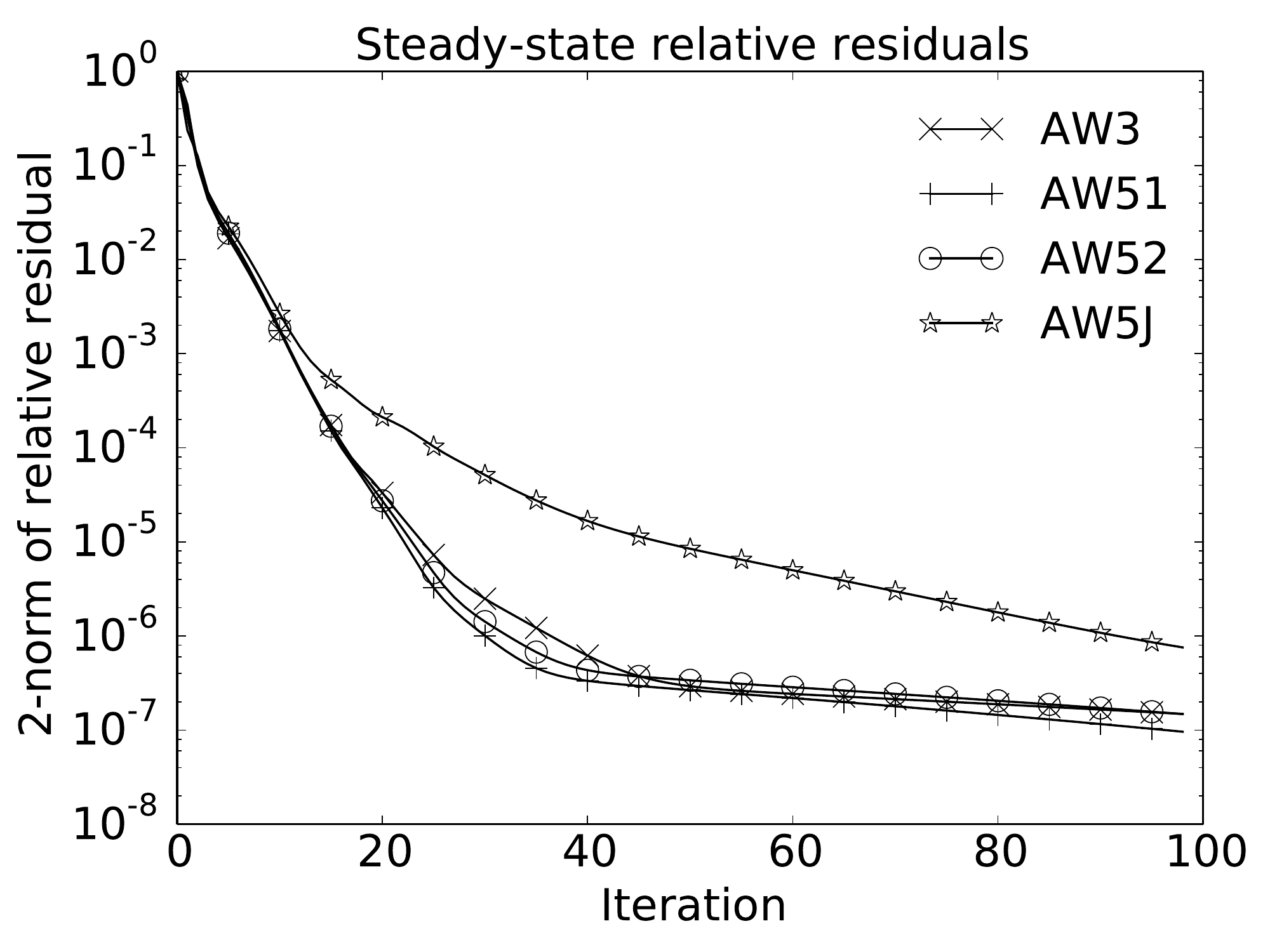}}
\subfigure[unsteady]{
\includegraphics[width=0.48\textwidth]{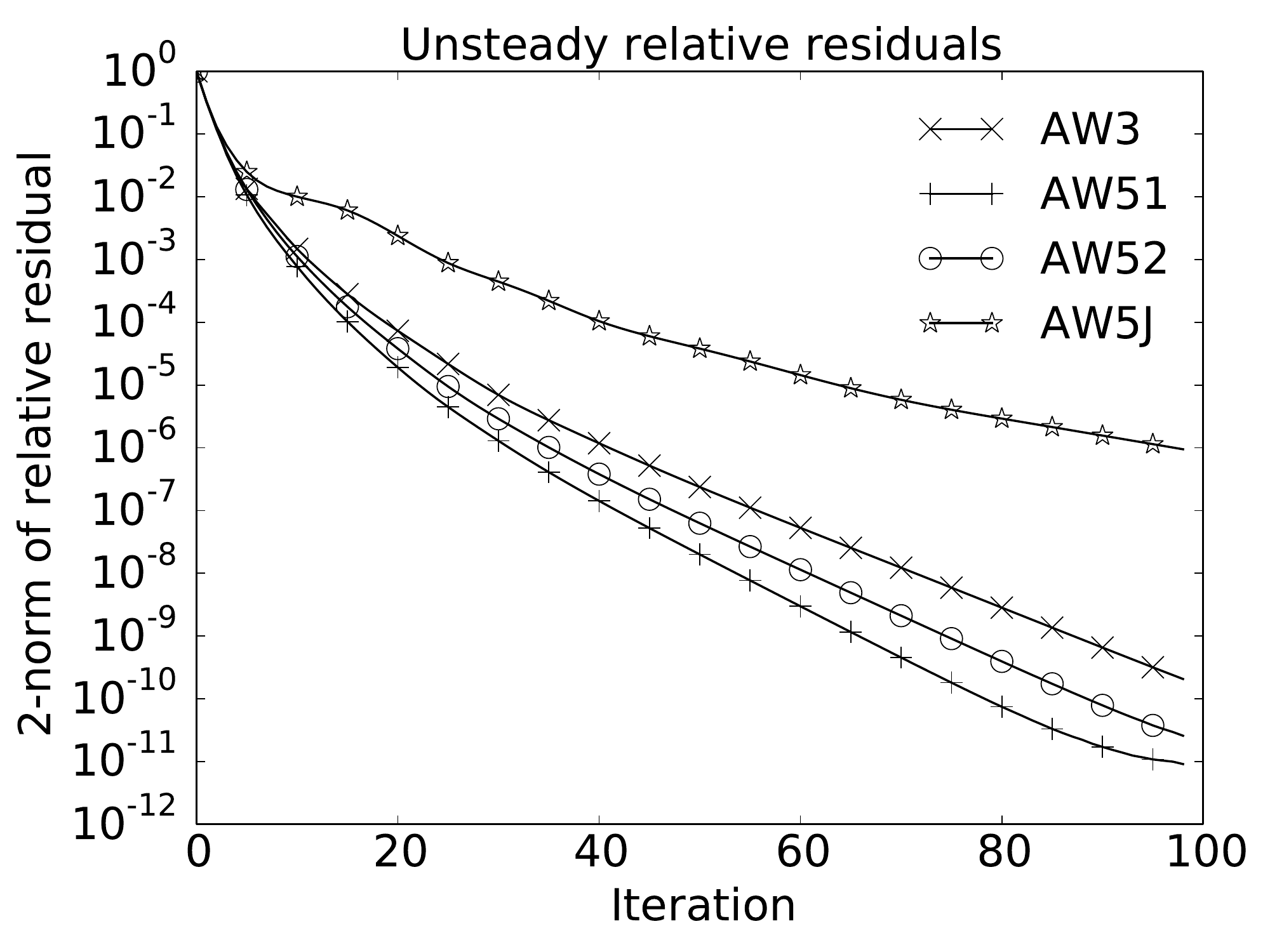}}
\caption{\label{conv_rae_all_d005}Convergence behavior for steady and unsteady flow around RAE 2822 airfoil for AW3, AW51 and AW52, $0^{\circ}$ angle of attack, 100 steady/unsteady iterations, $d$=0.05.}
\end{figure}

\begin{figure}[h]
\centering
\subfigure[steady]{
\includegraphics[width=0.48\textwidth]{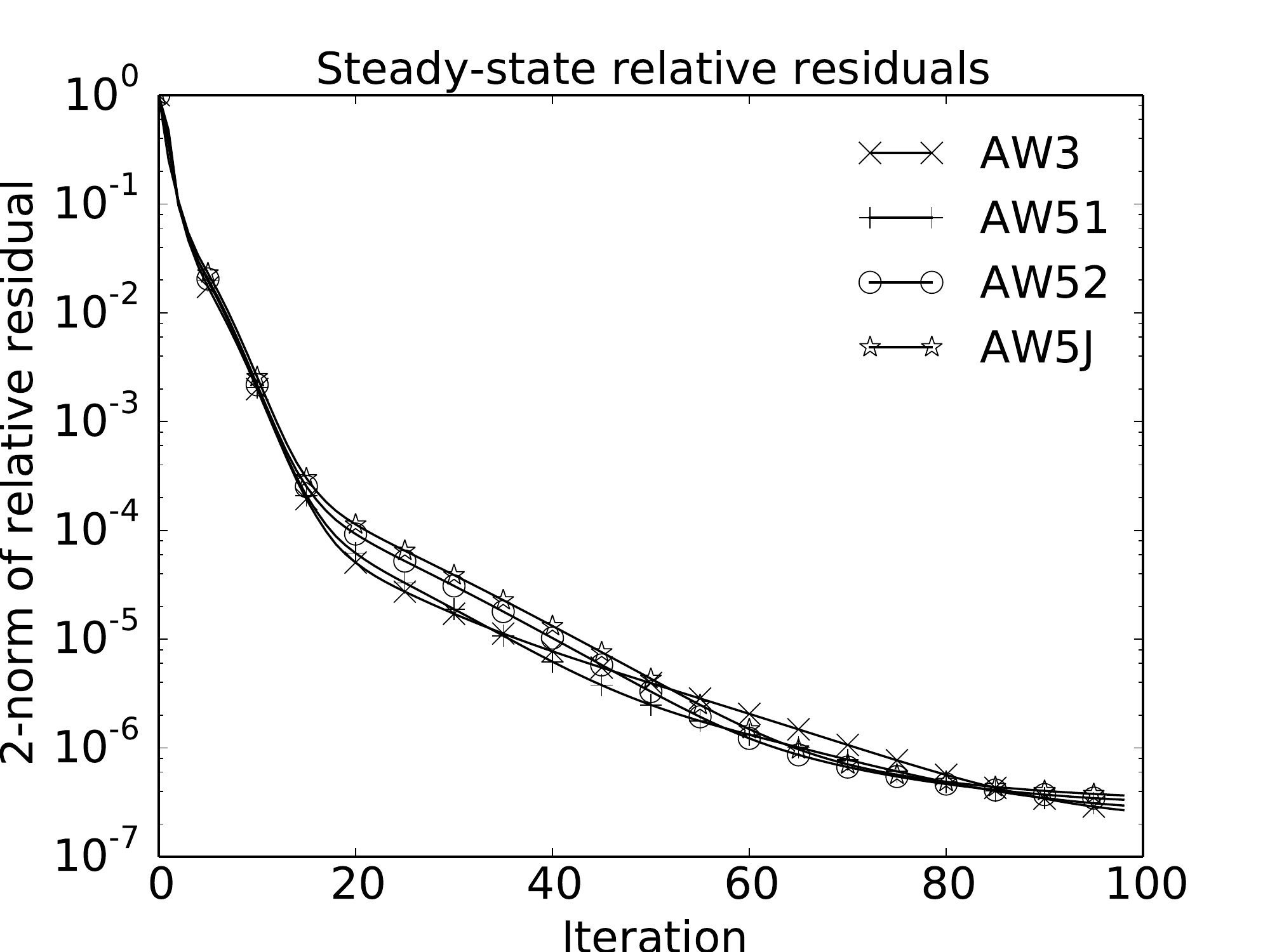}}
\subfigure[unsteady]{
\includegraphics[width=0.48\textwidth]{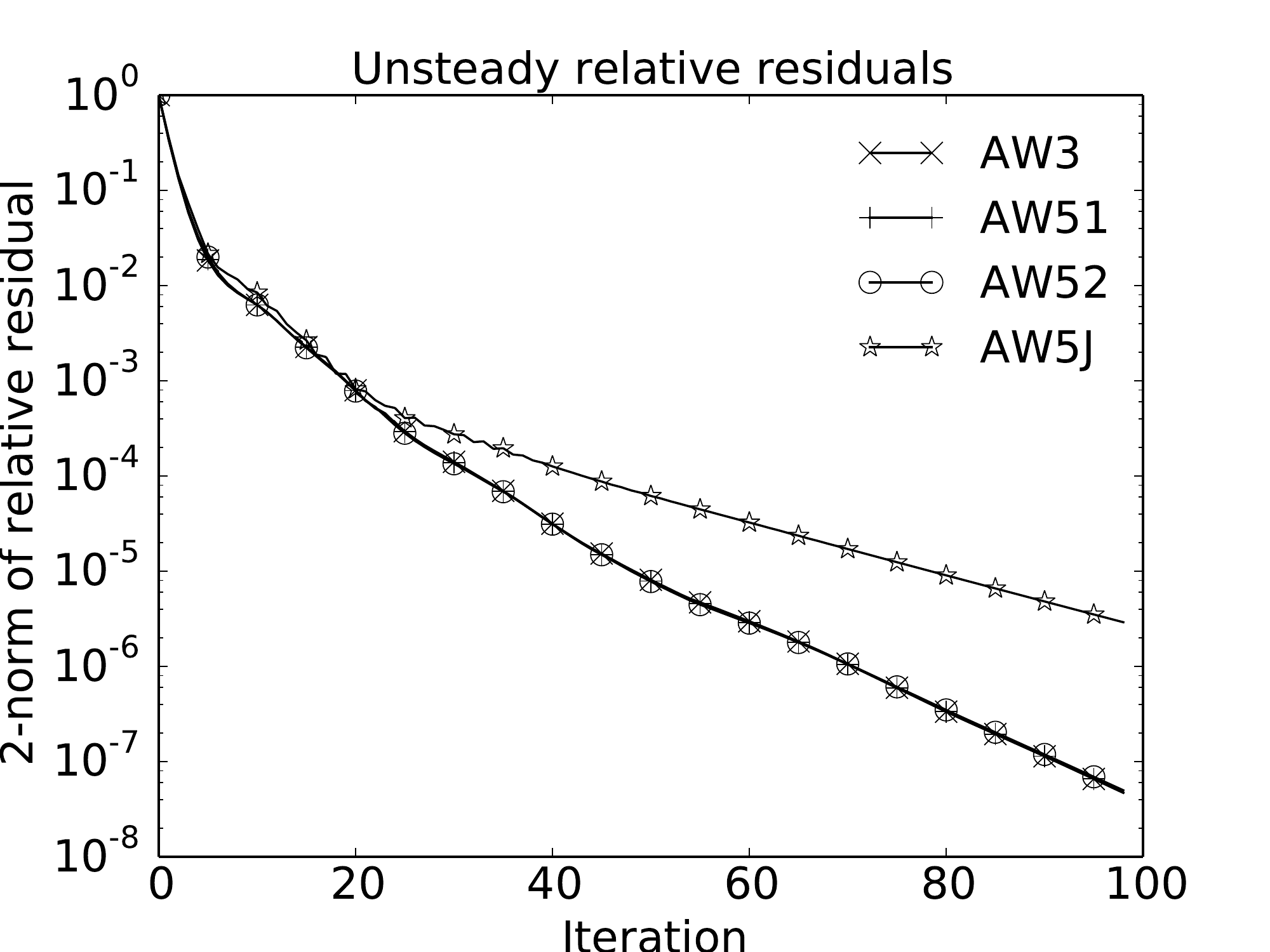}}
\caption{\label{conv_rae_all_d05}Convergence behavior for steady and unsteady flow around RAE 2822 airfoil for AW3, AW51 and AW52, $0^{\circ}$ angle of attack, 100 steady/unsteady iterations, $d$=0.5.}
\end{figure}

\begin{table}
\centering
\begin{tabular}{cc|cc|cc}
 & & \multicolumn{2}{|c|}{$d$=0.05, $c^*$=60/100} & \multicolumn{2}{|c}{$d$=0.5, $c^*$=10000} \\ \hline
       &        & CPU[s] & av. conv. rate & CPU[s] & av. conv. rate \\ \hline
Steady & AW3/60 & 8.6 & 0.8530 & 8.9 & 0.8581 \\
       & AW5J/20 & 12.7 & 0.8670 & 13.1 & 0.8609 \\
       & AW51/100 & 12.6 & 0.8492 & 13.1 & 0.8590 \\ 
       & AW52/80  & 12.7 & 0.8530 & 13.1 & 0.8601 \\ \hline
Unsteady & AW3 & 11.7  & 0.7972 & 12.0 & 0.8424 \\
BDF-2 & AW5J   & 17.5 & 0.8687 & 18.0 & 0.8786 \\ 
      & AW51   & 17.5 & 0.7732 & 18.1 & 0.8424 \\ 
      & AW52   & 17.9 & 0.7808 & 18.0 & 0.8429 \\ \hline
\end{tabular}
\caption{\label{uflo-rae}Performance of AW3, AW51 and AW52 for the RAE 2822 airfoil, $0^{\circ}$ angle of attack, 100 steady/unsteady iterations.}
\end{table}

The residual histories for the same tests, but for the RAE airfoil can be seen in Figures~\ref{conv_rae_all_d05} and \ref{conv_rae_all_d005}. Convergence rates and CPU times are summarized in table~\ref{uflo-rae}. The numbers after the scheme names are the values of $c^*$ used in the steady iterations. Again, faster convergence is obtained with $d$=0.05 for all schemes except AW5J.

As an immediate conclusion, it can be seen that the different schemes have similar convergence rates. Thus, AW3 performs best in terms of CPU times, since it is a three stage smoother, opposed to the five stage smoothers. With the fast mode, we get a convergence rate for the unsteady case of 0.77 for the NACA profile and 0.8 for the RAE profile. However, for the RAE profile, we have to reduce $c^*$ from 100 for 3 of the 4 schemes to prevent instability. With the convergence rate obtained, 20 iterations are sufficient for most applications, which is a matter of seconds. In the robust mode, the convergence rate goes down to 0.8 for the NACA profile and 0.84 for the RAE profile.

In the steady state case, there is a decline in convergence rate after 20 to 30 iterations. This explains why the convergence rates are significantly slower here. In the first phase, a convergence rate of about 0.7 is obtained and the norm of the residual is decreased by about $10^6$, which is completely sufficient for most applications.

\subsection{Mesh Independence}

To verify that the solvers' performance is mesh-independent, we run the pitching NACA 64010 airfoil with AW3 and $d$=0.5 on coarse ($256\times 32$), medium ($384\times 48$) and fine ($512\times 64$) meshes in robust mode. Table \ref{mesh-naca} shows the results. As can be seen, the convergence of the preconditioned W schemes is mesh-independent. With $d$ set to 0.05 the simulations on the coarse mesh diverged, which is an example where the robust mode is indeed more robust. 

\begin{table}
\centering
\begin{tabular}{c|cc}
Mesh & Steady & Unsteady \\ \hline
$256\times 32$ & 0.8698 & 0.8110 \\ 
$384\times 48$ & 0.8646 & 0.8144 \\ 
$512\times 64$ & 0.8616 & 0.8071\\ \hline
\end{tabular}
\caption{\label{mesh-naca}Convergence rate with AW3 smoothing for the pitching NACA 64A010 airfoil on different meshes, $c^*$=10000, $d=0.5$, pstep=36, $\alpha=0^{\circ}$.}
\end{table}

\subsection{Effect of flow angle}

In the Fourier analysis it was found that grid-aligned flow could be problematic. We therefore choose angles of attack $\alpha$ of 0, 1, 2 and 4 degrees for the steady state computation or the second time step in an unsteady computation. Table~\ref{angle-rae-naca} shows the convergence rates in fast mode ($d$=0.05). Essentially, it is unaffected by the angle of attack. However, for two cases, the iteration stalls after 30, resp. 65 iterations at relative residuals of $10^{-4}$ and $10^{-5}$, respectively. 

\begin{table}
\centering
\begin{tabular}{cc|cc}
       & Angle & NACA & RAE \\ \hline
Steady & 0     & 0.8623 & 0.8625 \\ 
       & 1     & 0.8500 & 0.8639 \\ 
       & 2     & 0.8597 & 0.8590 \\ 
       & 4     & 0.8530 & 0.8564 \\ \hline
Unsteady & 0   & 0.8249 & 0.8846 \\ 
       & 1   & 0.8239 & 0.9097 (65) \\ 
$\Delta t=0.486822$ & 2 & 0.8254 & 0.8307 \\ 
       & 4     & 0.9055 (30) & 0.8290 \\ \hline
\end{tabular}
\caption{\label{angle-rae-naca}Performance of AW3 for the NACA and RAE airfoil with varying $\alpha$, $c^*$=100, $d$=0.05, 100 steady/unsteady iterations.}
\end{table}

\section{Conclusions}

We considered preconditioned pseudo time iterations for agglomeration multigrid schemes for the steady and unsteady RANS equations. As a discretization, the JST scheme was used as a flux function in a finite volume method. We derived preconditioned additive W methods, as well as preconditioned additive explicit RK methods. Both are implemented in exactly the same way with the difference being in how the preconditioner is chosen, as well as the pseudo time step size. For the latter, the preconditioner has to approximate the Jacobian $\mathbf J$ and the pseudo time iteration has a finite stability region. In the additive W case, the preconditioner has to approximate $\mathbf I+\eta\Delta t^* \mathbf J$, whereby the pseudo time step size is possibly unbounded. However, we obtain an additional parameter $\eta$ which currently must be chosen empirically. As a preconditioner, we choose a flux vector splitting with a cutoff of small eigenvalues controlled by the free variable $d$. 

To compare the different methods, we used a discrete Fourier analysis of the linearized Euler equations. Numerical results show that AW3, AW51 and AW52 have similar convergence rates, meaning that AW3 performs best, since it uses two stages less. The free parameter $\eta$ can be chosen with relative freedom within a stable range ($0.5 \lessapprox \eta \lessapprox 0.9$) although the optimal value is dependent in some cases on the initial conditions, $d$ and the aspect ratio. Fixing $\eta=0.8$ is an acceptable simplification in the cases tested. The most significant parameter affecting stability and convergence is the eigenvalue cutoff coefficient $d$ in the numerical flux function. It was found that the W schemes were $A$-stable for $d \geq 0.5$ and had stability limits lower than preconditioned ARK schemes for $d < 0.5$. Thirdly, the pseudo CFL number $c^*$ was tuned for optimal performance. Different optimal values were obtained for different aspect ratios but as long as $c^*$ was within the stability limit, good convergence was achieved. This is useful since the aspect ratios in practical meshes vary considerably.

Simulations of pitching and plunging NACA 64A010 and RAE2822 airfoils in high Reynolds number flow at Mach 0.796 were performed using the 2D URANS code uflo103. The preconditioned ARK schemes were completely uncompetitive with convergence rates of around 0.999. The additive W schemes, on the other hand, achieved mesh-independent convergence rates of as low as 0.85 for the initial steady-state iteration and 0.77 for the unsteady iterations. Slightly different optimal values of $\eta$ and $c^*$ were found although the behaviour of the schemes was qualitatively similar to that predicted by the linear analysis. We emphasise two modes of operation for the W schemes: a fast mode, $d=0.05$, $\eta$=0.5 and $c^*$=100 and a robust mode, $d=0.5$, $\eta$=0.5 and $c^*$=10000. Unsteady convergence rates in the robust mode were higher than the fast mode but still competitive. Steady-state convergence rates for all tests stalled to varying degrees after around 20 iterations but the residuals had already fallen by 6 orders of magnitude - more than sufficient for most practical applications.

In summary, the new additive W schemes achieve excellent performance as smoothers in the agglomeration multigrid method applied to 2D URANS simulations of high Reynolds number transonic flows. The stiffness associated with very high aspect ratio grids is counteracted by highly tuned preconditioning. The underlying aim of this paper was to present a complete analysis of the reasons why such preconditioned iterative smoothers are effective, in order that their high performance can be replicated. We encountered two parameters that resisted analysis and had to be tuned empirically: $\eta$ and $d$. Nevertheless, this is considered a great improvement. Future work will look at these parameters in more detail. In addition, boundary conditions should have an influence on convergence speed. 

\section*{Acknowledgement}

We would like to thank Charlie Swanson for interesting discussions and sharing some code with us.

\end{document}